\documentclass[12pt]{article}
\usepackage{epsfig,wrapfig,epic}
\usepackage{amsmath}
\usepackage{amscd}
\usepackage{rotating}
\usepackage{wrapfig}

\usepackage{hyperref}
\usepackage{amsfonts}
\usepackage{mathrsfs}
\usepackage{stmaryrd}
\usepackage{shadow}
\usepackage{lineno}
\usepackage{color}

\usepackage{dsfont}
\usepackage{graphics}
\usepackage{graphicx}
\usepackage{amssymb}
\usepackage{amsmath}
\usepackage{url}

\DeclareMathAlphabet{\mathpzc}{OT1}{pzc}{m}{it}

\newcommand{\codim}{\mathpzc{codim}}
\newcommand{\dm}{\mathpzc{dim}}

\newcommand{\C}{\mathbb{C}}

\newcommand{\bdm}{\textbf{m}}
\newcommand{\bdn}{\textbf{n}}

\newcommand{\h}{{{\mbox{\scriptsize $\mathsf{H}$}}}}
\newcommand{\dist}{\mathpzc{dist}}
\newcommand{\dg}{\mathpzc{deg}}
\newcommand{\GCD}{\mathpzc{gcd}}
\newcommand{\spn}{\mathpzc{span}}
\newcommand{\lb}{\llbracket}
\newcommand{\rb}{\rrbracket}

\newcommand{\Jac}{\mathcal{J}}

\newtheorem{example}{Example}[section]
\newtheorem{Problem}{PROBLEM}[section]
\newtheorem{theorem}{Theorem}[section]

\newtheorem{cor}[theorem]{Corollary}
\newtheorem{lemma}[theorem]{Lemma}

\newtheorem{define}[theorem]{Definition}

\def\proof{\textsc{Proof.} }
\def\foorp{\hfill$\square$}

\title{The Numerical Factorization of Polynomials%
}

\author{ Wenyuan
Wu\thanks{Chongqing Institute of Green and Intelligent Technology,
Chinese Academy of Sciences. {\em Email:}
wuwenyuan@cigit.ac.cn. This work is partially supported by grant
NSFC 11471307. }\hspace{1cm} Zhonggang Zeng\thanks{Department of
Mathematics, Northeastern Illinois University, Chicago, Il 60625.
{\em Email:} zzeng@neiu.edu. Research supported in part by NSF under
Grant DMS-0715127.} }

\begin{document}

\input amssym.def
\maketitle

\vspace{-2mm}
\begin{abstract}
Polynomial factorization in conventional sense is an ill-posed problem
due to its discontinuity with respect to
coefficient perturbations, making it intractable for numerical computation using
empirical data.
~As a regularization, this paper formulates the notion of numerical
factorization based on the geometry of polynomial spaces and the stratification
of factorization manifolds.
~Furthermore, this paper establishes the existence, uniqueness, Lipschitz
continuity, condition number, and convergence of the numerical factorization
to the underlying exact factorization, leading to a robust and efficient
algorithm with a Matlab implementation capable of accurate polynomial
factorizations using floating point arithmetic even if the coefficients are
perturbed.
%
%
%
\end{abstract}



\vspace{-4mm}
\section{Introduction}

\vspace{-5mm}
Polynomial factorization is one of the fundamental algebraic operations in
theory and in applications.
~It is also an enduring research subject in the field of
computer algebra 
as well as a significant success of
symbolic computation (c.f. the survey \cite{kal03}).
~Factorization functionalities have been standard features of
computer algebra systems such as Maple and Mathematica with a
common assumption that the coefficients are represented exactly.
~Nonetheless, theoretical advancement and algorithmic development are still
in early stages in many cases.
~When a polynomial is approximately known with a limited
accuracy in coefficients, the very meaning of its factorization as we know it
becomes a question, as illustrated in the following example.
~More precisely, a well-posed notion of numerical factorization has not been 
established, leaving a gap in the foundation of its computation.

\vspace{-3mm}
\begin{example} \label{e:5fs}\em \footnotesize
~We illustrate the central question of this paper:
~{\em Assume the polynomial}
\begin{equation}\label{fxy}
 f ~=~ x^2 y+\mbox{\scriptsize 0.857143}\, x y^2
+\mbox{\scriptsize $0.833333$}\, x^2+\mbox{\scriptsize $1.38095$}\, x y
+\mbox{\scriptsize $0.571429$}\, y^2
+\mbox{\scriptsize $0.555556$}\, x+\mbox{\scriptsize $0.476190$}\, y
\end{equation}
{\em is given as the empirical data of a factorable polynomial
~$\tilde{f}$.
~Knowing that the data are imperfect with an error bound
~$\|f-\tilde{f}\| \le 10^{-5}$, ~what is the factorization of the underlying
polynomial ~$\tilde{f}\,$?}

The factorization of ~$f$ ~in conventional sense doesn't exist
while the underlying polynomial ~$\tilde{f}$ ~is factorable
but not known exactly.
~Intuitively, one can ask a more modest question:
{\em Is there a factorable polynomial near ~$f$ ~within the data error bound
~$10^{-5}$?}
~This latter question is similar to an open problem in \cite{kal99}
and the answer is ambiguous: ~The
polynomial ~$f$ ~is near many factorable
polynomials, as shown in Table~\ref{tb:fs}.

\begin{table}[h]
\begin{center}
\begin{tabular}{|l|c|}\hline
\multicolumn{1}{|c|}{factorable polynomial near ~$f$ }
& distance \\ \hline
\raisebox{-2mm}{%
$\tilde{f}_{~} ~=~~  (x+\mbox{\scriptsize $\frac{2}{3}$})
(y+\mbox{\scriptsize $\frac{5}{6}$}) (x+\mbox{\scriptsize $\frac{6}{7}$}\,y)$}
\hfill \raisebox{-2mm}{\scriptsize $\longleftarrow$ underlying polynomial
~~~~~}
&
\raisebox{-2mm}{\scriptsize $2.53\times 10^{-6}$}  \\
$\hat{f}_{~}~=~$
$\mbox{\tiny $0.9999994$}\,
(x + \mbox{\tiny $0.6666667$})
(y + \mbox{\tiny $0.8333333$})
(x + \mbox{\tiny $0.8571429$}\,y)$
~\scriptsize $\longleftarrow$ the numerical factorization
&
\mbox{\scriptsize $2.11\times 10^{-6}$}  \\
$f_1~=~$
$\mbox{\tiny $1.0000002$}\,
(y + \mbox{\tiny $0.8333327$})
(x^2 + \mbox{\tiny $0.6666663$}\,x +
\mbox{\tiny $0.5714287$}\,y +
\mbox{\tiny $0.8571420$}\,xy)$
&
\mbox{\scriptsize $1.75\times 10^{-6}$}  \\
$f_2~=~$
$\mbox{\tiny $1.0000002$}\,
(x + \mbox{\tiny $0.8571425$}\,y)
(xy + \mbox{\tiny $0.8333321$}\,x +
\mbox{\tiny $0.6666661$}\,y +
\mbox{\tiny $0.5555555$})
$ &
\mbox{\scriptsize $1.60\times 10^{-6}$}  \\
$f_3~=~ $
$\mbox{\tiny $0.9999997$}\,
(x + \mbox{\tiny $0.66666728$})
(xy + \mbox{\tiny $0.83333331$}\,x +
\mbox{\tiny $0.7142836$}\,y + \mbox{\tiny $0.8571432$}\,y^2)$
&
\mbox{\scriptsize $8.70\times 10^{-7}$}
\\
\hfill
\raisebox{1mm}{
\begin{rotate}{-15} \scriptsize $\longleftarrow$ \end{rotate}}
~~~~\scriptsize the nearest factorable polynomial & \\
\hline
\end{tabular}
\end{center}
   \vspace{-5mm}
\caption{\small The polynomial ~$f$ ~in (\ref{fxy}) is near many factorable
polynomials with various distances.} \label{tb:fs}
\end{table}

The {\em numerical factorization} of ~$f$ ~within the error
tolerance ~$10^{-5}$, ~as we shall define in \S\ref{s:nf}, is the exact
factorization of ~$\hat{f}$ ~in Table~\ref{tb:fs} and
accurately approximates the factorization of
~$\tilde{f}$ ~from which ~$f$ ~is constructed by rounding up digits.
~The nearest polynomial to the data ~$f$, ~however,
is not ~$\hat{f}$ ~but ~$f_3$ ~whose factorization does {\em not} resemble
that of ~$\tilde{f}$.
~In fact, the factorable polynomial with the smallest distance to the
data is almost certain to have an incorrect factorization structure
by the Factorization Manifold Embedding Theorem in \S\ref{s:str} whenever the
underlying polynomial ~$\tilde{f}$ ~has more than two factors.
\foorp
\end{example}

\vspace{-3mm}
As shown in this example, conventional factorization is a so-called
{\em ill-posed problem} for numerical computation since the factorization is
discontinuous with respect to data perturbations.
~Consequently, fundamental questions arise such as if, under what conditions,
by computing which factorization and to what accuracy
we can recover the factorization from empirical data.
%
%
%
~In this paper, we establish the geometry of the polynomial (topological) spaces
in Factorization Manifold Theorem and Factorization Manifold
Embedding Theorem.
~Based on the geometry we rigorously formulate the notion of the
{\em numerical factorization}.
~We prove the so-defined numerical factorization eliminates the ill-posedness
of the conventional factorization and accurately approximates the intended
exact factorization (Numerical Factorization Theorem) with a finite
sensitivity measure that is conveniently attainable (Numerical Factorization
Sensitivity Theorem).
~As a result, the intractable ill-posed factorization problem in numerical
computation is completely regularized as a well-posed numerical factorization
problem that approximates the intended factorization with an accuracy in the
same order of the data precision.

Our results can be narrated as follows.
~The collection of polynomials possessing a nontrivial factorization structure 
is a complex analytic manifold of a positive codimension and every such
manifold is embedded in the closures of certain manifolds of
lower codimensions.
~This dimension deficit provides a singularity measurement
of polynomials on the manifold and fully explains the ill-posedness their
factorizations:
~An infinitesimal perturbation reduces the singularity
and pushes a polynomial away from its native manifold into the
open dense subset of polynomials with a trivial factorization structure,
making the exact factorization on the empirical data meaningless.
~Based on the geometric analysis, we formulate the notion of the
numerical factorization as the exact factorization of the polynomial
on the nearby factorization manifold of the highest singularity
having the smallest distance to the data.
~Under the assumption that the data error
is small, the original factorization can be recovered
accurately by the numerical factorization of the data polynomial within a
proper error bound even if it is perturbed.
~From the Tubular Neighborhood Theorem in differential geometry,
the numerical factorization is a well posed problem as it uniquely
exists, is Lipschitz continuous and approximates the exact factorization
of the underlying polynomial the data represent.
~The accuracy of the recovered factorization is in the same
order of the data accuracy since the factorization is
Lipschitz continuous on that manifold.
~Moreover, the conventional factorization becomes a special case of the
numerical factorization within a small error tolerance.
~The analysis of numerical factorization leads to a
two-staged computing strategy for the numerical factorization:
~Identifying the factorization manifold by a squarefree factorization and a
proper reducibility test, followed by the Gauss-Newton iteration \cite{DS02,zz05}
for minimizing the distance to the factorization manifold.

This paper attempts to bridge differential geometry, computer algebra and 
numerical analysis.
~As an effective analytical tool that still appears to be underused, 
geometry has led to many penetrating insights in numerical analysis
(e.g. \cite{BC13,Kahan72}) and effective algorithms such as homotopy methods
based on Sard's Theorem and Theorem of Bertini (e.g. \cite{bertini,Li03}).
~Polynomial factorization problem has been studied from geometric perspective
such as in \cite{CGHKW01,CGKW02,GalHoe}.
~This paper broadens the geometric analysis into a numerical computation
of a basic problem in computer algebra by establishing the stratified complex 
analytic manifolds of factorization and their tubular neighborhood.
~In a seminal technical report \cite{Kahan72}, Kahan is the first to discover
the hidden continuity on manifolds for generally discontinous solutions of
ill-posed algebraic problems.
~Recent works such as \cite{zz05,Zeng09b} made progress along this directions.
~This work provides a complete regularization of a typical ill-posed algebraic 
problem in numerical polynomial factorization by establishing its existence,
uniqueness, Lipschitz continuity, convergence and condition number.
~Regularizations to this extent should now be expected for other ill-posed 
algebraic problems that share a similar geometry.

For exact polynomial factorization, many
effective methods
have been developed over the past several decades.
~Those algorithms and complexity analyses have been studied
extensively.
~The work of Sasaki Suzuki, Kolar and Sasaki \cite{54} introduces the
techniques
of extended Hensel construction and the trace recombination
that lead to factorization algorithms
such as van Hoeij's trace recombination \cite{23} for univariate polynomial
factorization of integer coefficients.
~The first polynomial-time factorization algorithms
is given by Lenstra, Lenstra and Lovasz \cite{44} for univariate polynomial
factorization, and by
Kaltofen and Von zur Gathen \cite{20,31} for multivariate polynomials.
~Rigorous proofs are also provided in these works
on the probabilities and the complexities.
~At present, the algorithm having the lowest complexity for
exact bivariate polynomial factorization appears to be due to
Lecerf \cite{43}.

Many authors made pioneer contributions to the numerical
factorization problem of multivariate polynomials, such as
pseudofactors by Huang, Stetter, Wu and Zhi \cite{stetter},
the numerical reducibility tests by Galligo and Watt \cite{GalWat}
and by Kaltofen and May \cite{kal-may},
computing zero sum relations by Sasaki \cite{Sasaki01}, interpolating
the irreducible factors as curves by Corless, Giesbrecht, Van Heij,
Kotsireas and Watt \cite{CGHKW01},
and by
Corless, Galligo, Kotsireas and Watt \cite{CGKW02}.
~Finding a nearby factorable polynomial as proposed in
\cite{CGHKW01,GalWat,GalHoe,GKMYZ04,kal99,kal03,kal-may}
has played an indispensable role in
the advancement of numerical polynomial factorization, even though such a
backward accuracy alone is insufficient in numerical factorizations
as illustrated in Example~\ref{e:5fs}.
In \cite{SVW04}, Sommese, Verschelde and Wampler developed a homotopy
continuation method along with monodromy grouping,
and Verschelde released and has maintained the first numerical factorization
software as part of the PHC package \cite{PHC} for solving polynomial systems.
%
~A breakthrough due to Ruppert's differential forms \cite{rup99} led
to a novel hybrid factorization algorithm \cite{Gao03}
by Gao, and the development of a numerical factorization
algorithm by Gao, Kaltofen, May, Yang and Zhi in \cite{GKMYZ04,KMYZ08}.
~Based on the formulation and analysis of this paper, we developed an
algorithm that shares a root similar to \cite{Gao03,GKMYZ04,KMYZ08} along with
several new developments as well as a Matlab implementation.

The results of this paper is not limited to
multivariate polynomials.
~The numerical factorization theory and computational strategy extend to
the univariate polynomial factorization, which is also known as polynomial
root-finding where a recent major development enables accurate
computation of multiple roots without extending the hardware precision
even if the coefficients are perturbed \cite{zz05}.
%
%
~This paper provides a unified framework for
the numerical factorization including the univariate
factorization as a special case.

\vspace{-5mm}
\section{Preliminaries}

\vspace{-5mm}
We consider polynomials in variables ~$x_1,\,\ldots\,,x_{\ell}$
~with coefficients in the field ~$\C$ ~of complex numbers.
~The ring of these polynomials is commonly denoted by
~$\mathbb{C}[x_1,\,\ldots\,,x_{\ell}]$.
~The ~$\ell${\bf -tuple degree} of a polynomial ~$f$ ~is defined as a vector
~$\dg(f) ~=~ \big(\dg_{x_1}(f),\,\ldots\, , \dg_{x_{\ell}}(f)\big)$
~where ~$\dg_{x_j}(f)$ ~is the degree of ~$f$ ~in ~$x_j$.
%
%
%
%
%
%
%
~For any ~$\ell$-tuple degree ~$\mathbf{n}$, ~denote
\begin{eqnarray*}
  \mathbb{P}^{\textbf{n}} &~:=~&
\big\{ p \in \mathbb{C}[x_1,\,\ldots\,,x_{\ell}]
~\big|~ \dg(p) \leq \textbf{n} \big\} \\
 \mathscr{P}^{\textbf{n}} &~:=~&
\big\{ p \in \mathbb{C}[x_1,\,\ldots\,,x_{\ell}]
~\big|~ \dg(p) = \textbf{n} \big\}.
\end{eqnarray*}
Here ~$\mathbb{P}^{\textbf{n}}$ ~is a vector space whose dimension is
denoted by ~$\langle\textbf{n}\rangle$. 
~Inequality between ~$\ell$-tuple degrees are componentwise.
%
%
~With a monomial basis in lexicographical order,
a polynomial
~$f = f_1 x^{\bdm_1}+f_2 x^{\bdm_2}+\cdots+f_{\langle \bdn \rangle}
x^{\bdm_{\langle \bdn \rangle}}$ ~in
~$\mathbb{P}^\bdn$ ~corresponds to a unique {\bf coefficient vector}
denoted by
~$\lb f \rb := (f_1,\ldots,f_{\langle \bdn \rangle}) ~\in
\C^{\langle \bdn \rangle}$,
~such as ~$f= 3\,x_1^2x_2 - 4\,x_1x_2 + 5\,x_1 + 6 \in
\mathbb{P}^{(2,1)}$ ~corresponding to
~$\lb f \rb  ~=~
(3,0,-4,5,0,6) ~\in~ \C^6$.
~A subset ~$\Omega \subset \mathbb{P}^\mathbf{m}$ ~corresponds to the
subset ~$\lb \Omega \rb \,=\,\big\{ \lb p \rb \in
\C^{\langle \mathbf m \rangle}
~\big|~ p \in \Omega \big\}$ ~in ~$\C^{\langle \mathbf m \rangle}$.
~Here ~$\C^n$ ~is the vector space of ~$n$-dimensional vectors of complex
numbers.
~All vectors in this paper are ordered arrays denoted by
boldface lowercase letters or in the form of ~$\lb \cdot \rb$.
~The Euclidean norm in ~$\mathbb{C}^{\langle \mathbf{n} \rangle}$ ~induces
the polynomial norm as
~$\|f\|\,:=\, \big\| \lb f \rb \big\|_2$, ~making ~$\mathbb{P}^\mathbf{n}$
~a topological metric space.
%
%
%
%

There are no differences between factoring a polynomial
and factoring its nonzero constant multiple.
~We say ~$p$ ~and ~$q$ ~are
~{\bf equivalent}, denoted by ~$p\sim q$,
~if ~$p = \alpha q$ ~for ~$\alpha\in\C\setminus\{0\}$.
~A metric is needed in the quotient space
~$\C[x_1,\ldots,x_\ell]/\sim$ ~but not seen in the literature.
~We propose a scaling-invariant distance between
polynomials ~$p$ ~and ~$q$ ~as the sine of the principal
angle between the subspaces ~$\spn\{ \lb p \rb \}$ ~and ~$\spn\{ \lb q \rb\}$,
~denoted by
\begin{equation} \label{sine}
  \sin(p,q) ~~:=~~ \left\{ \begin{array}{ccl}
0 & & \mbox{if}~~~ p = q = 0 \\
1 & & \mbox{if}~~~ p = 0, ~q \ne 0 ~\mbox{or}~ p \ne 0, ~q = 0 \\
\left\| \frac{p}{\|p\|} -
\frac{\lb q \rb \cdot \lb p \rb}{\|q\|\,\|p\|}
\frac{q}{\|q\|} \right\| & & \mbox{if}~~~p \ne 0, ~q \ne 0.
\end{array} \right.
\end{equation}
Here the ~``$\cdot$'' ~denotes the standard vector dot product.
~Let ~$P_f$ ~be the projection mappings to
~$\spn\{\lb f \rb\}$ ~for any polynomial ~$f$.
~It is known that ~$\sin(p,q) \equiv \|P_p -
P_q\|_2$ ~(c.f. \cite{stew}), ~and is thus
a distance in the quotient space ~$\C[x_1,\ldots,x_l]/\sim$.

A polynomial ~$f$ ~is
{\bf factorable} if there exist nonconstant polynomials ~$g$ ~and ~$h$
~such that ~$f=g\, h$, ~otherwise it is
{\bf irreducible}.
~We say ~$\alpha f_1\,f_2\,\cdots\,f_k$ ~is a
{\bf factorization} of ~$f$ ~if ~$\alpha \in \C \setminus \{0\}$,
~$\dg(f_j) \ne \textbf{0}$ ~for ~$j=1,\ldots,k$, ~and
~$\alpha f_1 \cdots f_k \,\sim\,f$.
~Here we abuse the notation ~$\alpha f_1\cdots f_k$ ~as it represents
either the polynomial product or the factorization
that consists of factors ~$\alpha$, \,$f_1,\, \ldots,\, f_k$
~depending on the context.
~We say two factorizations ~$\alpha\, f_1 \, f_2 \, \cdots\, f_k$
~and ~$\beta\, g_1 \, g_2 \, \cdots\, g_m$ ~are equivalent, denoted by
~$\alpha\, f_1 \, f_2 \, \cdots\, f_k \,\sim\,
\beta\, g_1 \, g_2 \, \cdots\, g_m$,
~if ~$m=k$ ~and
there is a permutation ~$\{\sigma_1,\ldots,\sigma_k\}$ ~of ~$\{1,\ldots,k\}$
~such that ~$f_j \sim g_{\sigma_j}$ ~for ~$j=1,\ldots,k$.
~If ~$f_1$, $f_2$, $\cdots$, $f_k$ ~are all irreducible,
then ~$\alpha f_1\,f_2\,\cdots\,f_k$ ~is an {\bf irreducible
factorization}.
%
~The irreducible factorization of a polynomial is unique as an equivalence class.

A factorization ~$\gamma g_1 \cdots g_m$ ~is regarded
as an approximate factorization of ~$f$ ~if the {\bf
backward error} ~$\sin(f,\gamma g_1 \cdots g_m)$ ~is small enough
and acceptable in the underlying application.
~The {\bf forward error} of the factorization ~$\gamma g_1 \cdots g_m$
~is the difference between the factors ~$g_1,\,\ldots,\,g_m$ ~and their
counterparts in ~$f = \alpha\, f_1 \cdots f_k$ ~via a proper metric that
is needed but not properly established in the literature.
~Here we extend the distance measurement ~$\sin(\cdot,\cdot)$ ~to the
{\bf distance between two factorizations} as
\begin{equation}
\dist\big(\alpha f_1 \cdots f_k, ~
\gamma\, g_1 \cdots g_m\big) ~~:=~~
\left\{ \begin{array}{cl}
\mbox{\raisebox{1mm}{$1$}} & ~~\mbox{\raisebox{1mm}{if ~$m \ne k$}}  \\
{\displaystyle \min_{(\sigma_1,\ldots,\sigma_k) \in \Sigma}}
\Big\{~{\displaystyle \max_{1\le j \le k}}
\big\{\sin (f_j,g_{\sigma_j})\big\} \Big\} &
~~\mbox{otherwise}
\end{array} \right.
\label{ferr}
\end{equation}
where ~$\Sigma$ ~is the collection of all permutations
~$(\sigma_1,\ldots, \sigma_k)$ ~of ~$(1,\ldots,k)$.
~Clearly, two factorizations are equivalent if and only if their distance
is zero.

A polynomial is {\bf squarefree} if its irreducible factorization
consists of pairwise coprime factors.
~A {\bf squarefree factorization} ~$\alpha f_1^{k_1} \cdots f_r^{k_r}$
~consists of squarefree polynomials ~$f_1, \ldots, f_r$ ~as components that
are pairwise coprime but may or may not be irreducible.
%
%
~Again, we use the notation
~$\alpha f_1^{k_1}f_2^{k_2} \cdots f_r^{k_r}$
~to represent either the polynomial that equals to the result of
the polynomial multiplication or the factorization consists of the
factors ~$\alpha$, ~$f_1,\ldots,f_1, f_2,\ldots,f_2, \ldots, f_r,\ldots,f_r$
~where each ~$f_j$ ~repeats ~$k_j$ ~times for ~$j=1,\ldots,r$.
~If ~$\alpha f_1^{k_1}\cdots f_r^{k_r}$ ~is an irreducible squarefree
factorization of a polynomial ~$f$ ~with degree ~$\mathbf{m}$,
~we shall use
~$\mathfrak{M} ~=~ \textbf{m}_{1}^{k_1}\cdots\textbf{m}_{r}^{k_r}$
~to denote the {\bf factorization structure}, or simply the {\bf structure}
of ~$f$, ~where
~$\textbf{m}_{j} = \dg(f_j) \ne \mathbf{0}$,
~$k_j \ge 1$ ~for ~$j=1,\ldots,r$
~and ~$k_1\mathbf{m_1}+\cdots +k_r \textbf{m}_r = \textbf{m} = \dg(f)$.
~We shall also say such an ~$\mathfrak{M}$ ~is one of the factorization
structures of the degree ~$\mathbf{m}$ ~and denote
~$\dg(\mathfrak{M}) = \mathbf{m}$.
~Any permutation
of ~$\mathbf{m}_1^{k_1}, \,\ldots\,, \mathbf{m}_r^{k_r}$ ~in
~$\mathfrak{M} ~=~ \textbf{m}_{1}^{k_1}\cdots\textbf{m}_{r}^{k_r}$ ~is
considered the same structure.
~There are two cases for a factorization structure ~$\mathfrak{M}$
~to be called {\bf trivial} when ~$\mathfrak{M}$ ~is the
factorization structure of either an irreducible polynomial or a
univariate polynomial with no multiple roots.
%
%
%
%
~A factorization structure is {\bf nontrivial} if it is not trivial.

\vspace{-5mm}
\section{Factorization Manifolds} \label{s:geo}

\vspace{-5mm}
The factorization of a polynomial ~$f$ ~is an equivalence
class in which a specific representative ~$\alpha f_1^{k_1}\cdots f_r^{k_r}$
~can be extracted using a set of auxiliary equations ~$\textbf{b}_1 \cdot
\lb f_{1} \rb \,=\, \cdots \,=\, \textbf{b}_r \cdot \lb f_r \rb \,=\, 1$
~where ~$\textbf{b}_1, \ldots, \textbf{b}_r$ ~are unit vectors of proper
dimensions.
~We call such vectors ~$\textbf{b}_1,\ldots,\textbf{b}_r$ ~the
{\bf scaling vectors}.
%
%
~Scaling vectors can be chosen randomly.
~A more natural choice during computation is the
normalized initial approximation of ~$\lb f_i \rb$ ~so that
~$\mathbf{b}_i \cdot \lb f_i \rb \,\approx\, \|f_i\|^2 = 1$ ~for
~$i=1,\ldots,r$.
~For any factorizations ~$\gamma p_1^{k_1}\cdots p_r^{k_r}$ ~and ~$\mu
q_1^{k_1}\cdots q_r^{k_r}$ ~scaled by equations
~$\mathbf{b}_i \cdot \lb p_i \rb = \mathbf{b}_i \cdot \lb q_i \rb  =
1$ ~for ~$i=1,\ldots,r$, ~it is clear that ~$\|p_1\|, \ldots, \|p_r\|,
\|q_1\|, \ldots, \|q_r\| \ge 1$ ~since the scaling vectors are of unit norms,
and the following lemma applies.

\vspace{-3mm}
\begin{lemma} \label{l:dis}
~Let ~$\gamma p_1^{k_1}\cdots p_r^{k_r}$ ~and ~$\mu q_1^{k_1}\cdots
q_r^{k_r}$ ~be two factorizations with
~$\|p_i\|, \|q_i\| \ge 1$ ~for
~$i=1,\ldots,r$.
~Then
\begin{equation} \label{disin}
\dist\big(\gamma p_1^{k_1}\cdots p_r^{k_r},\,\mu q_1^{k_1}\cdots q_r^{k_r}
\big) ~~\le~~ \max_{1\le i \le r} \big\|p_i - q_i\big\|.
\end{equation}
\end{lemma}

\vspace{-3mm}
\proof
~It is straightforward to verify that ~$\sin(p_i,q_i) \le \|p_i - q_i\|$
~whenever ~$\|p_i\|, \|q_i\| \ge 1$ ~for ~$i=1,\ldots,r$.
~Thus (\ref{disin}) holds.
\foorp

Suppose ~$f$ ~possesses an irreducible squarefree factorization
~$\alpha f_1^{k_1} \cdots f_r^{k_r}$ ~and ~$\dg(f_i) = \textbf{m}_{i}$
for ~$i=1,\ldots,r$.
~The factorization structure
~$\mathfrak{M}$ ~equals to ~$\textbf{m}_{1}^{k_1}\cdots\textbf{m}_{r}^{k_r}$.
~All the polynomials sharing this factorization structure form
a subset
\begin{eqnarray*} \mathscr{F}^\mathfrak{M}  ~~:=~~
\big\{ f \in \mathscr{P}^\mathbf{m}  &~\big|~ &
f = \alpha g_1^{k_1}\cdots g_r^{k_r} ~~\mbox{where}~\,
\alpha \in \C,  ~g_j \in \mathscr{P}^ {\textbf{m}_{j}},
~j=1,\ldots,r  \\
& & ~~\mbox{are irreducible and pairwise coprime}
\big\}
\end{eqnarray*}
of ~$\mathbb{P}^{\textbf{m}}$ ~where
~$\textbf{m} = \dg(\mathfrak{M})$.
~For almost all unit scaling vectors
~$\mathbf{b}_i \in \C^{\langle \mathbf{m}_i \rangle}$ ~for
~$i=1,\ldots,r$,
~a polynomial ~$f\in \mathscr{F}^\mathfrak{M}$ ~possesses irreducible
factors ~$\alpha, f_1,\ldots,f_r$ ~such that
~$\alpha f_1^{k_1}\cdots f_r^{k_r}  \,=\, f$ ~and
~$\textbf{b}_1 \cdot \lb f_1 \rb = \cdots =
\textbf{b}_r \cdot \lb f_r \rb~= 1$
~so that the array ~$\big(\gamma,\lb p_1 \rb, \ldots, \lb p_r \rb\big)
= \big(\alpha, \lb f_1\rb,\ldots,\lb f_r\rb\big)$
~is a solution to the equation
\begin{equation} \label{phieq}
\phi\big(\gamma, \lb p_1\rb,\ldots,\lb p_r\rb\big) ~~=~~
\big(\lb f\rb, 1, \ldots, 1\big)
\end{equation}
for
~$\gamma \in \C$, ~$p_j \in \mathbb{P}^{\mathbf{m}_j}$, ~$j=1,\ldots,r$,
~where the mapping ~$\phi$ ~is defined by
\begin{equation} \label{eq:diffeo2}
\begin{array}{rcl}
\phi ~:~ \C \times \C^{\langle \mathbf{m}_1 \rangle} \times \cdots \times
\C^{\langle \mathbf{m}_r \rangle} &~\longrightarrow~&
\C^{\langle \mathbf{m} \rangle} \times \C \times \cdots \times \C \\
\big(\gamma, \lb p_1 \rb,\ldots,\lb p_r \rb\big) &~\longmapsto~&
 \big( \lb \gamma p_1^{k_1}\cdots p_r^{k_r} \rb,\,
    \textbf{b}_1 \cdot \lb p_1 \rb, \, \ldots , \,
    \textbf{b}_r \cdot \lb p_r \rb
\big).
\end{array}
\end{equation}%
Let ~$q_i = (\frac{k_i}{p_i})\, \alpha p_1^{k_1} 
\cdots p_{r}^{k_r}$.
~Then the Jacobian of ~$\phi$ ~can be written as
\begin{eqnarray}
\lefteqn{\Jac(\alpha, \lb p_1\rb,\ldots,\lb p_r \rb) ~~=} \nonumber \\
& & \left[
\begin{array}{ccccc}
\mathpzc{column}( \lb p_1^{k_1}\cdots p_r^{k_r} \rb) & C_{\textbf{m}_{1}}(q_1) &
C_{\textbf{m}_{2}}(q_2) &
\cdots & C_{\textbf{m}_{r}}(q_r) \\
& \mathpzc{column}(\textbf{b}_1)^\h &  &  &  \\
&  & \mathpzc{column}(\textbf{b}_2)^\h  &  & \\
&  &    &  \ddots &   \\
&  &    &  & \mathpzc{column}(\textbf{b}_r)^\h \\
                    \end{array}
                  \right]\label{eq:Jac_diffeo2}
\end{eqnarray}
where ~$\mathpzc{column}{(\cdot)}$ ~represents
the column block generated by a 
vector ~$(\cdot)$, ~the notation ~$(\cdot)^\h$ ~denotes
the Hermitian transpose of the matrix ~$(\cdot)$,
~and ~$C_{\textbf{m}_{i}}(q_i)$
~is the convolution matrix \cite{zz05} ~associated with ~$q_i$ ~so that
~$C_{\textbf{m}_{i}}(q_i) \cdot \lb h \rb ~=~ \lb q_i h \rb$
~holds for any ~$h \in \mathbb{P}^{\mathbf{m}_i}$, ~$i=1,\ldots,r$.
~We need several lemmas for establishing the main theorems of the paper.

\vspace{-3mm}
\begin{lemma}\label{lem:injective}
~For ~$\alpha \in \C\setminus \{0\}$,
~$\mathbf{b}_j \in \C^{\langle \mathbf{m}_j\rangle}$ ~and
~$p_j \in \mathscr{P}^{\mathbf{m}_j}$ ~with
~$\mathbf{b}_j \cdot \lb p_j \rb \ne 0$ ~for ~$j=1,\ldots,r$,
~the Jacobian in {\em (\ref{eq:Jac_diffeo2})} is injective if
and only if  ~$p_1, \ldots, p_r$ ~are pairwise coprime.
\end{lemma}

\vspace{-3mm}
\proof
~Assume ~$p_1,\ldots,p_r$ ~are pairwise coprime and the
matrix-vector multiplication
%
\begin{equation}\label{eq:Jac_inj2}
\Jac(\alpha, \lb p_1\rb,\ldots,\lb p_r \rb) \cdot
   \left( - a, \, \lb v_1 \rb, \, \ldots, \, \lb v_r \rb \,
                                \right) ~~=~~ \textbf{0}.
\end{equation}
Then
~$\textbf{b}_1 \cdot \lb v_1 \rb = \cdots = \textbf{b}_r \cdot \lb v_r \rb = 0$
~as well as
~$\sum_{i=1}^r q_i v_i \,=\, a\prod_{i=1}^r p_i^{k_i}$
~that lead to ~$
\sum_{i=1}^r k_i \alpha p_1\cdots p_{i-1} v_i p_{i+1}\cdots p_r \,=\,
a p_1 p_2 \cdots p_r$.
Thus ~$k_1 \alpha v_1 p_2\cdots p_r$ $=$ $a p_1p_2\cdots p_r-$ $
\sum_{i=2}^r  k_i \alpha p_1\cdots p_{i-1} v_i p_{i+1}\cdots p_r$
~that contains the factor ~$p_1$.
~Because ~$\GCD(p_1,p_j) = 1$ ~for ~$j=2,\,\ldots\,,r$, ~there is a
polynomial ~$s$ ~such that ~$v_1 = s p_1$.
~The degree ~$\dg(v_1) \le \dg(p_1)$
~leads to ~$s$ ~being a constant.
~Since ~$\textbf{b}_1 \cdot \lb p_1 \rb \neq 0$, ~$0 = \textbf{b}_1 \cdot
\lb v_1 \rb = s\, \textbf{b}_1 \cdot \lb p_1 \rb$, ~hence ~$s = 0$.
~Consequently ~$v_1 = 0$.
~Similarly we can prove that ~$v_i = 0$ ~for ~$i=2,\ldots,r$.
~Substituting ~$v_1 = \cdots = v_r = 0$ ~into (\ref{eq:Jac_inj2}),
we have ~$a
p_1^{k_1} \cdots p_r^{k_r} = 0$ ~and thus ~$a=0$.
~Therefore, the Jacobian is injective.
~Conversely, to prove that the injectiveness of the Jacobian in
(\ref{eq:Jac_diffeo2}) implies ~$p_1, \ldots , p_r$
~are pairwise coprime,
assume there are some ~$i \ne j$ ~such that ~$\GCD(p_i,p_j) \ne 1$.
~Then we shall prove that the Jacobian must be rank-deficient.
~Without loss of generality, we can assume ~$p_1
= e\,s$ ~and ~$p_2 = e\,t$ ~for some polynomials ~$e$, ~$s$ ~and ~$t$
~where ~$e = \GCD(p_1,p_2)$ ~is nonconstant.
~Then there are three possible cases.
~As case one, if ~$\textbf{b}_1 \cdot \lb s \rb = 0$ ~and
~$\textbf{b}_2 \cdot \lb t \rb = 0$, ~then it is easy to show that
~$\big(0,\,\frac{1}{k_1}\lb s \rb,\,\frac{-1}{k_2}\lb t \rb,\,\textbf{0},
\,\ldots\,,\textbf{0}\big)$
~is a
nonzero solution to (\ref{eq:Jac_inj2}).
~As case two, if ~$\textbf{b}_1 \cdot \lb s \rb = c \neq 0$ ~and
~$\textbf{b}_2 \cdot \lb t \rb = 0$, ~then we can consider ~$w =
\frac{c}{\beta_1}\, e\,s - s$, ~where ~$\beta_1 = \textbf{b}_1 \cdot
\lb p_1 \rb$.
~It is straightforward to verify that ~$\textbf{b}_1 \cdot \lb w \rb =
\frac{c}{\beta_1}\, \textbf{b}_1 \cdot \lb p_1 \rb - \textbf{b}_1 \cdot
\lb s \rb = c - c = 0$.
~Since ~$e  =  \GCD(p_1,p_2)$ ~which is nontrivial, we have ~$\dg(ces)
> \dg(s)$ ~and consequently ~$w \neq 0$.
~Thus ~$\big(\frac{-c}{\beta_1}\alpha, \,\frac{1}{k_1}\lb w \rb,
\,\frac{1}{k_2}\lb t \rb,
\,\textbf{0},\,\ldots\,,\textbf{0}\big)$
~is a nonzero solution of (\ref{eq:Jac_inj2}).
~For the third case where ~$\textbf{b}_1 \cdot \lb s \rb = c \neq 0$
~and ~$\textbf{b}_2 \cdot \lb t \rb = d \neq 0$, ~let ~$v_1 =
\frac{c}{\beta_1}\; e\,s - s$ ~and ~$v_2 = -\frac{d}{\beta_2}\; e\,t +
t$ ~where ~$\beta_2 = \textbf{b}_2 \cdot \lb p_2 \rb$.
~Then  ~$\big(\frac{d\alpha}{\beta_2}-\frac{c\alpha}{\beta_1} ,
\,\frac{1}{k_1}\lb v_1 \rb, \,\frac{1}{k_2}
\lb v_2 \rb,\,\textbf{0},\,\ldots\,,\textbf{0}\big)$
~is a nonzero solution of (\ref{eq:Jac_inj2}).
~Therefore, the Jacobian is a rank-deficient matrix. \foorp

\vspace{-3mm}
\begin{lemma} \label{l:pcv}
~Let ~$\mathfrak{M} = \mathbf{m}_1^{k_1}\cdots\mathbf{m}_r^{k_r}$
~be a factorization structure of degree ~$\mathbf{m}$
~and assume a sequence
~$\{p_j\}_{j=1}^\infty \subset \mathscr{F}^\mathfrak{M}$
~converges to ~$q\in \mathscr{P}^\mathbf{m}$.
~Then there is a subsequence of ~$\{p_j\}_{j=1}^\infty$
~whose irreducible factorizations
converge to a factorization ~$\alpha q_1^{k_1}\cdots q_r^{k_r}$
~of ~$q$ ~with ~$\dg(q_i) = \mathbf{m}_i$ ~for ~$i=1,\ldots,r$.
~Further assume ~$q \in \mathscr{F}^\mathfrak{M}$.
~Then the irreducible factorizations of
~$\{p_j\}_{j=1}^\infty$
~converge to the irreducible factorization of ~$q$.
\end{lemma}

\vspace{-3mm}
\proof
~Let ~$\mathfrak{S}_i = \big\{ \lb  f \rb ~\big|~ f \in
\mathbb{P}^{\mathbf{m_i}}, ~\|f\|=1 \big\}$ ~and
~$p_j = \alpha_j p_{j1}^{k_1}\cdots p_{jr}^{k_r}$
~be an irreducible factorization
of ~$p_j$, ~where $p_{ji} \in \mathfrak{S}_i$ ~for ~$i \in \{1,\,\ldots\,,r\}$
~and ~$j=1,2,\,\ldots$.
~Denote ~$\textbf{P}_j = (\alpha_j, \lb p_{j1}
\rb,\,\ldots\,, \lb p_{jr} \rb) \in \C \times
\mathfrak{S}_1 \times \cdots \times \mathfrak{S}_r$.
~There is a subsequence ~$\{j_1,j_2,\ldots\}$ ~of ~$\{1,2,\ldots\}$
~such that
~$\displaystyle \lim_{\sigma\rightarrow\infty} \|p_{j_\sigma i} - q_i\| = 0$
~for ~$i \in \{1,\,\ldots\, ,r\}$ ~since ~$\mathfrak{S}_i$'s are compact.
~As a result,  the subsequence ~$\big\{\alpha_{j_\sigma}\big\}$
~converges to certain ~$\alpha\in\C$.
%
Namely, the subsequence
~$\{\textbf{P}_{j_\sigma}\}_{\sigma=1}^{\infty}$ ~converges to a point
~$(\alpha, \lb q_1 \rb,\ldots, \lb q_r \rb)$
~such that ~$q=\alpha q_1^{k_1}\cdots q_r^{k_r}$
~and ~$\dg(q_i) \leq \mathbf{m}_i$ ~for ~$i=1,\ldots,r$.
~From ~$\dg(q) = \mathbf{m}$ ~we have ~$\dg(q_i) =\mathbf{m}_i$ ~for
~$i=1,\ldots,r$.
~By Lemma~\ref{l:dis}, the irreducible factorizations of
~$p_{j_\sigma}$ ~for ~$\sigma=1,2,\ldots$ ~converge to the factorization
~$\alpha q_1^{k_1}\cdots q_r^{k_r}$
~since $\dist{(\alpha_{j_\sigma} p_{j_\sigma\,1}^{k_1}\cdots
p_{j_\sigma\,r}^{k_r}, ~\alpha q_1^{k_1}\cdots q_r^{k_r})}$ $\leq
\max_{i} \|p_{j_\sigma\,i}-q_i\| \longrightarrow 0$
~when ~$\sigma \rightarrow \infty$.
~Moreover,  if ~$q \in \mathscr{F}^\mathfrak{M}$, ~then
~$\alpha q_1^{k_1}\cdots q_r^{k_r}$ ~is an irreducible squarefree
factorization of ~$q$ ~by the uniqueness of factorizations.
~Furthermore, the irreducible squarefree factorizations of
the whole sequence ~$\{p_j\}$ ~must converge to the factorization
~$\alpha  q_1^{k_1}\cdots q_r^{k_r}$
~since otherwise there would be a ~$\delta > 0$ ~and
a subsequence of ~$\{\textbf{P}_j\}_{j=1}^{\infty}$ ~converging to
~$(\hat{\alpha}, \hat{q}_1,\ldots, \hat{q}_r)$, ~with
~$q = \hat{\alpha} \hat{q}_1^{k_1}\cdots \hat{q}_r^{k_r}$
~and ~$\dist \big(\hat{\alpha} \hat{q}_1^{k_1}\cdots \hat{q}_r^{k_r},
\alpha  q_1^{k_1}\cdots q_r^{k_r} \big) \ge \delta$,
~contradicting the uniqueness of the factorization of ~$q$.
 \foorp

Lemma~\ref{l:pcv} directly leads to the following corollaries.

\vspace{-3mm}
\begin{cor}\label{c:phi}
~Let ~$f$ ~be a polynomial with a factorization structure
~$\mathfrak{M} = \mathbf{m}_1^{k_1}\cdots \mathbf{m}_r^{k_r}$
~of degree ~$\mathbf{m}$ ~and an irreducible squarefree
factorization ~$\alpha f_1^{k_1} \cdots f_r^{k_r}$ ~satisfying
~$\|f_1\|$ $=\cdots=$ $\|f_r\|=1$.
~For any ~$\epsilon >0$, ~there is a neighborhood ~$\Omega_{f}$ ~of ~$f$ ~in
~$\mathbb{P}^\mathbf{m}$ ~such that every ~$g \in \Omega_f
\cap \mathscr{F}^\mathfrak{M}$ ~corresponds to a unique
~$(\beta, g_1, \ldots,g_r)$
~with ~$g = \beta g_1^{k_1} \cdots g_r^{k_r}$,
~$\lb f_1 \rb \cdot \lb g_1 \rb = \cdots =
\lb f_r \rb \cdot \lb g_r \rb = 1$ ~and
~$\sqrt{|\alpha-\beta|^2+\|f_1-g_1\|^2+\cdots+\|f_r-g_r\|^2}
\,<\, \epsilon$.
\end{cor}

\vspace{-3mm}
\proof
~For any ~$\delta > 0$, ~Lemma~\ref{l:pcv} implies that there is a
neighborhood ~$\Omega_{f,\delta}$ ~of ~$f$ ~in ~$\mathbb{P}^\mathbf{m}$
~such that the irreducible squarefree factorization
~$\beta g_1^{k_1} \cdots g_r^{k_r}$ ~of
every $g \in \Omega_{f,\delta} \cap \mathscr{F}^\mathfrak{M}$
~satisfies ~$\dist\big(\beta g_1^{k_1} \cdots g_r^{k_r},
~\alpha f_1^{k_1} \cdots f_r^{k_r}\big) < \delta$.
~We can assume ~$\delta < \frac{1}{2} \min_{i\ne j}\{ \sin(f_i,f_j)\}$
~and ~$\max_j\,\{\sin(f_j,g_j)\} < \delta$.
~Since ~$\sin(f_i,g_j) \ge \sin(f_i,f_j)-\sin(f_j,g_j)$ $> \delta$
~whenever ~$i\ne j$, ~no other permutation of ~$g_1,\ldots,g_r$ ~satisfies
~$\max_j\{\sin(f_j,g_j)\} < \delta$.
~Further assume ~$g_1,\ldots,g_r$ ~are the unique representatives in
their respective equivalence classes satisfying ~$\lb f_j \rb
\cdot \lb g_j \rb = 1$ ~for ~$j=1,\ldots,r$.
~Then ~$\lb f_j \rb \cdot \lb f_j - g_j \rb = 0$,
~$\|f_j-g_j\| = \|g_j\|\sin(f_j,g_j)$ ~and ~$\|g_j\|^2 = \|f_j\|^2 +
\|f_j-g_j\|^2$ $=1+\|g_j\|^2\sin^2(f_j,g_j)$, ~leading to
~$\|f_j-g_j\| = \frac{\sin(f_j,g_j)}{\sqrt{1-\sin^2(f_j,g_j)}} <
\frac{\delta}{\sqrt{1-\delta^2}}$ ~for ~$j=1,\ldots,r$.
~Therefore, for any ~$\epsilon > 0$,
~the assertion holds when ~$\delta$ ~is small.
\foorp

\vspace{-3mm}
\begin{cor}\label{cor:irred}
~Polynomials of degree ~$\mathbf{m}$ ~with a trivial
factorization structure form an open subset of
~$\mathbb{P}^{\mathbf{m}}$.
\end{cor}

\vspace{-3mm}
\proof
~For a univariate degree ~$\mathbf{m}$, ~the assertion follows from
the continuity of polynomial roots with respect to the coefficients.
~Assume ~$\mathbf{m}$ ~is multivariate and the assertion does not hold.
~Then there is an irreducible polynomial ~$f$ ~of degree ~$\mathbf{m}$
~and a sequence of factorable polynomials ~$\{p_j\}_{j=1}^{\infty}$
~approaching ~$f$.
~Because there are finitely many factorization
structures in ~$\mathbb{P}^{\mathbf{m}}$,
~there exists a nontrivial factorization structure
~$\mathfrak{M}$ ~and
a subsequence ~$\{p_{j_\sigma}\}_{\sigma=1}^\infty$
~in ~$\mathscr{F}^\mathfrak{M}$.
~By Lemma \ref{l:pcv}, the irreducible factorizations of this
subsequence converge to a nontrivial
factorization of ~$f$, ~contradicting the
irreducibility of ~$f$.
\foorp

We can now establish the following Factorization Manifold Theorem.
~A subset ~$S$ ~in the topological space ~$\mathbb{P}^\mathbf{m}$ ~is a
{\bf complex analytic manifold}
of dimension ~$k$ ~in ~$\mathbb{P}^\mathbf{m}$
~if, for every
~$p \in S$,
~there exists an open subset ~$\Omega$ ~of ~$\mathbb{P}^\mathbf{m}$
~containing ~$p$ ~and a biholomorphic mapping from
~$\lb S\cap \Omega\rb \subset \C^{\langle \mathbf{m}\rangle}$ ~onto
an open subset of ~$\C^k$.
~The \textbf{codimension}, namely the dimension deficit, of
~$S$ ~is denoted by ~$\codim(S) :=
\dm(\mathbb{P}^{\mathbf{m}}) - \dm(S) =
\langle\mathbf{m}\rangle - k$.
~The Factorization Manifold Theorem is at core of the geometry on the 
polynomial factorization.
~This result and the proof are fundamental but not seen in the literature.

\vspace{-3mm}
\begin{theorem}[Factorization Manifold Theorem]
\label{th:codim} ~Let ~$\mathfrak{M} =
\mathbf{m}_1^{k_1}\cdots\mathbf{m}_r^{k_r}$ ~be a factorization
structure of degree ~$\mathbf{m}$.
~Then ~$\mathscr{F}^\mathfrak{M}$ ~is a complex analytic manifold
in ~$\mathbb{P}^{\mathbf{m}}$ ~and
\begin{equation} \label{eq:codim}
 \codim (\mathscr{F}^\mathfrak{M}) ~~=~~
\langle \mathbf{m} \rangle
 - \Big(\langle \textbf{m}_{1} \rangle + \cdots +
\langle \textbf{m}_{r} \rangle +1 - r\Big).
\end{equation}
\end{theorem}

\vspace{-3mm}
\proof
~Let ~$f  \in \mathscr{F}^{\mathfrak{M}}$ ~with a irreducible
squarefree factorization ~$\alpha f_1^{k_1}\cdots f_r^{k_r}$ ~where
~$\dg(f_j) = \mathbf{m}_j$ ~and ~$\|f_1\|=\cdots=\|f_r\|=1$.
~Setting
~$\mathbf{b}_j = \lb f_j \rb$ ~for ~$j=1,\ldots,r$ ~in
(\ref{eq:diffeo2}) yields a holomorphic mapping ~$\phi$ ~from
~$\C^k$ ~to ~$\C^{\langle\mathbf{m}\rangle+r}$ ~with
~$k=1+\langle\mathbf{m}_1\rangle+\cdots+\langle\mathbf{m}_r\rangle$
~and ~$\phi\big(\alpha, \lb f_1 \rb, \ldots, \lb f_r \rb\big) =
\big(\lb f\rb, 1, \ldots, 1\big)$.
~By Corollary~\ref{cor:irred}, there is a neighborhood
~$\Delta$ ~of ~$\big(\alpha,\lb f_1\rb, \ldots, \lb f_r\rb\big)$
~in ~$\C\times \C^{\langle \mathbf{m}_1\rangle} \times \cdots \times
\C^{\langle \mathbf{m}_r\rangle}$ ~and
every ~$\big(\tilde{\alpha}, \lb \tilde{f}_1\rb ,\ldots,\lb \tilde{f}_r\rb\big)
\in \Delta$
~forms an irreducible squarefree factorization
~$\tilde{\alpha} \tilde{f}_1^{k_1} \cdots \tilde{f}_r^{k_r}
\in \mathscr{F}^{\mathfrak{M}}$.
%
~By Lemma~\ref{lem:injective}, the Jacobian
of ~$\phi$ ~is of full
rank ~$k$ ~at ~$(\alpha, \lb f_1 \rb,\ldots,\lb f_r\rb)$.
~As a result, the Inverse Mapping Theorem 
ensures
that certain ~$k$ ~components of ~$\phi$ ~form a biholomorphic
mapping ~$\check{\phi}$ ~from an open neighborhood ~$\Sigma$ ~of
~$\big(\alpha,\lb f_1 \rb, \ldots, \lb f_r \rb\big)$ ~in ~$\C^k$ ~to
an open subset ~$\Pi$ ~of ~$\C^k$.
~We can assume ~$\Sigma \subset \Delta $.
%
~This ~$\check{\phi}$ ~must contain the last ~$r$ ~components of
~$\phi$ ~since ~$\check\phi$ ~would not be injective without those
scaling constraints.
~Without loss of generality, we assume ~$\check\phi$ ~consists
of the last ~$k$ ~components of ~$\phi$ ~and we split
~$\phi(\mathbf{x})$ ~into
~$\hat\phi(\mathbf{x}) ~=~ \mathbf{u}$, ~and
~$\check\phi(\mathbf{x}) ~=~ \big( \mathbf{v}, \, \mathbf{w} \big)$
where ~$\mathbf{u} \in \C^{\langle\mathbf{m}\rangle+r-k}$,
~$\mathbf{v} \in \C^{k-r}$ ~and ~$\mathbf{w} \in \C^{r}$.
~Let ~$\check\Pi ~=~ \left\{ \mathbf{v} \in \C^{k-r} ~\big|~
\big( \mathbf{v}, 1,\ldots,1\big)
\in \Pi \right\}$ ~which is open in ~$\C^{k-r}$.
~Then the mapping ~$\mu(\mathbf{v}) = \big(
\hat\phi\circ\check\phi^{-1}(\mathbf{v},1,\ldots,1), \,
\mathbf{v} \big)$
~defined from ~$\check\Pi$ ~to ~$\C^{\langle\mathbf{m}\rangle}$
~is holomorphic,
and ~$\mu(\tilde\Pi) \subset \mathscr{F}^{\mathfrak{M}}$
~since ~$\tilde\Pi \times \{(1,\ldots,1)\} \subset \check\Pi$
~and ~$\Sigma \subset \Delta$.
%
%
~Furthermore, define ~$\psi ~:~  \C^{\langle\mathbf{m}\rangle+r-k} \times
\check\Pi \longrightarrow \check\Pi$ ~as the projection
~$\psi(\mathbf{u},\mathbf{v}) = \mathbf{v}$.
%
%
~By Corollary~\ref{c:phi},
there is an open neighborhood ~$\Omega \subset \C^{\langle \mathbf{m} \rangle
-k+r} \times \check{\Pi}$ ~of ~$\lb f \rb$ ~in
~$\C^{\langle \mathbf{m}\rangle}$ ~such that every
~$\lb p\rb \in \Omega \cap \lb \mathscr{F}^{\mathfrak{M}}\rb$
~corresponds to a unique ~$(\gamma, \lb p_1 \rb, \ldots, \lb p_r \rb)
\in \Sigma$ ~with
~$\phi(\gamma,\lb p_1 \rb, \ldots, \lb p_r \rb ) \,=\,
\big( \lb p \rb,\, 1,\, \ldots,\, 1 \big)$, ~namely
~$\mu\circ \psi(\lb p \rb) = \lb p \rb$.
~Define ~$\tilde{\Pi} = \psi(\Omega \cap \lb \mathscr{F}^{\mathfrak{M}}\rb)$.
~We have ~$\mu(\tilde\Pi) \supset \Omega \cap \lb
\mathscr{F}^{\mathfrak{M}}\rb)$.
~Then for every ~$\mathbf{v} \in \tilde{\Pi}$, ~there is a ~$\mathbf{u}
\in \C^{\langle \mathbf{m}\rangle+r-k}$ ~such that ~$(\mathbf{u,v}) \in
\Omega \cap \lb \mathscr{F}^{\mathfrak{M}}\rb$ ~corresponds to a
unique ~$(\gamma,\lb p_1 \rb, \ldots, \lb p_r \rb ) \in \Sigma$ ~with
~$\phi(\gamma,\lb p_1 \rb, \ldots, \lb p_r \rb ) =
(\mathbf{u,v},1,\ldots,1)$, ~implying ~$\mathbf{u} =
\hat\phi\circ \check\phi^{-1}(\mathbf{v},1,\ldots,1)$ ~and ~thus
~$\mu(\mathbf{v}) = (\mathbf{u,v})$.
~Namely ~$\mu(\tilde\Pi) \subset
\Omega \cap \lb \mathscr{F}^{\mathfrak{M}}\rb$ ~and thus
~$\mu(\tilde\Pi) =
\Omega \cap \lb \mathscr{F}^{\mathfrak{M}}\rb$.
~Since ~$\Omega$ ~is open and ~$\mu$ ~is continuous, hence
~$\mu^{-1}(\Omega) = \mu^{-1}(\Omega \cap \mathscr{F}^{\mathfrak{M}}) =
\tilde\Pi
$ ~is open in ~$\C^{k-r}$.
~Therefore, ~$\psi$ ~is biholomorphic from
$\Omega \cap \lb \mathscr{F}^{\mathfrak{M}}\rb$ ~onto
~$\tilde\Pi$ ~with the inverse ~$\mu$.
~Namely ~$\mathscr{F}^{\mathfrak{M}}$ ~is a complex analytic manifold of
dimension ~$k-r$, ~and (\ref{eq:codim}) follows.
\foorp

We shall refer to  ~$\mathscr{F}^\mathfrak{M}$ ~as the {\bf factorization
manifold} associated with the factorization structure ~$\mathfrak{M}$.
~Its dimension deficit indicates
how ill-posed the factorization is for polynomials on the
manifold.
~For a polynomial ~$p$ ~of degree ~$\mathbf{m}$,
~we say the {\bf singularity} ~of ~$p$ ~and its
factorization structure ~$\mathfrak{M}$ ~is ~$k$ ~if
~$\mathscr{F}^\mathfrak{M}$ ~is of
codimension ~$k$ ~in ~$\mathbb{P}^\mathbf{m}$.
~A polynomial is {\bf singular} in terms of factorization if its
singularity is positive, or {\bf nonsingular} otherwise.

\vspace{-3mm}
\begin{cor} \label{c:od}
~A polynomial is singular if and only its factorization structure
is nontrivial, and nonsingular polynomials of degree ~$\mathbf{m}$
~form an open dense subset of ~$\mathbb{P}^\mathbf{m}$.
\end{cor}

\vspace{-3mm}
\proof ~For both type of trivial factorization structures,
the corresponding factorization
manifold has a singularity zero
from (\ref{eq:codim}) by a straightforward verification.
~To prove ~$\codim(\mathscr{F}^\mathfrak{M}) > 0$ ~for
any nontrivial structure
~$\mathfrak{M}$, ~it suffices to show that for any ~degrees
~$\hat{\mathbf{n}} \ne \mathbf{0}$ ~and ~$\check{\mathbf{n}} \ne \mathbf{0}$
~such that ~$\hat{\mathbf{n}}+\check{\mathbf{n}}$ ~is a non-univariate degree,
we have
\begin{equation} \label{dineq}
\langle\hat{\mathbf{n}}+\check{\mathbf{n}}\rangle ~~>~~
\langle\hat{\mathbf{n}}\rangle + \langle\check{\mathbf{n}}\rangle -1.
\end{equation}
In fact, if is straightforward to verify (\ref{dineq}) for ~$\ell = 2$,
~namely we have the inequality
~$(\hat{n}_1+\check{n}_1+1)(\hat{n}_2+\check{n}_2+1) >
(\hat{n}_1+1)(\hat{n}_2+1) + (\check{n}_1+1)(\check{n}_2+1) -1$
~if ~$\hat{n}_1+\check{n}_1 > 0$ ~and ~$\hat{n}_2+\check{n}_2 > 0$,
~and the inequality (\ref{dineq}) for any positive integer ~$\ell\ge 2$
~follows an induction.

Let ~$\mathfrak{M}$ ~be trivial.
~$\mathscr{F}^{\mathfrak{M}}$ ~is open in ~$\mathbb{P}^\mathbf{m}$
~by Corollary \ref{cor:irred}.
~It is dense in ~$\mathbb{P}^\mathbf{m}$
~since it equals
~$\mathbb{P}^\mathbf{m}$ ~minus finitely many singular factorization manifolds
of lower dimensions.
\foorp

Corollary \ref{c:od} provide an ultimate explanation
why polynomial factorization is an ill-posed problem:
~Any polynomial ~$p$ ~having a nontrivial factorization
is singular in terms of factorization.
~Almost all perturbations ~$\Delta p$ ~results in
~$\tilde{p} = p + \Delta p$ ~that is pushed off the native manifold
into the open dense subset of nonsingular polynomials,
altering the factorization to a trivial one.
~This discontinuity makes the conventional
factorization ill-posed and intractable in numerical computation.
~When the factorization structure is preserved, however,
the irreducible factorization is Lipschitz continuous as asserted in the
following corollary.
~It is this continuity that makes numerical factorization possible.

\vspace{-3mm}
\begin{cor}[Factorization Continuity Theorem] \label{c:iflip}
~The irreducible factorization is locally Lipschitz continuous
on a factorization manifold:
~For any polynomial ~$f \in \mathscr{F}^\mathfrak{M}$
~with an irreducible squarefree factorization
~$\alpha f_1^{k_1} \cdots f_r^{k_r}$,
~there are constants ~$\delta, \eta > 0$ ~such that, for every polynomial
~$g\in \mathscr{F}^\mathfrak{M}$ ~satisfying ~$\|f-g\| < \delta$,
~the irreducible squarefree factorization
~$\beta g_1^{k_1} \cdots g_r^{k_r}$ ~of ~$g$ ~satisfies
~$ \dist\big(\alpha f_1^{k_1} \cdots f_r^{k_r},
~\beta g_1^{k_1} \cdots g_r^{k_r}\big)
~\le~ \eta\, \|f-g\|$.
\end{cor}

\vspace{-3mm}
\proof
~We can assume ~$\|f_1\|=\cdots=\|f_r\|=1$ ~and thus define the mapping
~$\phi$ ~in (\ref{eq:diffeo2}) with ~$\mathbf{b}_j = \lb f_j \rb$ ~for
~$j=1,\ldots,r$.
~Using the notations in the proof of Theorem~\ref{th:codim}, the
mapping ~$\zeta(\lb g\rb) =
\check{\phi}^{-1}\big(\psi(\lb g \rb ),1,\ldots,1\big)$
~is holomorphic and thus
Lipschitz continuous for any
~$g \in \mathscr{F}^\mathfrak{M}$ ~near ~$f$.
~Thus the assertion of this corollary follows from Lemma~\ref{l:dis}.
\foorp


\vspace{-5mm}
\section{Geometry of Factorization Manifolds} \label{s:str}

\vspace{-5mm}
Factorization manifolds form a topologically
stratified space ~$\mathbb{P}^\mathbf{m}$
~in which every singular factorization manifold is embedded in manifolds of
lower singularities as we shall elaborate in detail.
~There are two embedding operations on a factorization structure:
~The {\bf degree combining} operation is adding two \,$\ell$-tuple
degrees of the
same multiplicity while keeping other components of the factorization
structure unchanged:
\begin{equation} \label{dgcmb}
\cdots\mathbf{n}_i^{k_i} \cdots \mathbf{n}_j^{k_j}\ldots
~~\longrightarrow~~
\cdots (\mathbf{n}_i+\mathbf{n}_j)^{k} \cdots
~~~~\mbox{where}~~~ k_i ~=~ k_j ~=~ k.
\end{equation}
The {\bf multiplicity splitting} operation decomposes a component of a
factorization structure into two as follows:
\begin{equation} \label{msplt}
\cdots \mathbf{n}_i^{k_i} \cdots
~~\longrightarrow~~
\cdots \mathbf{n}_i^{\hat{k}_i} \,\mathbf{n}_i^{\tilde{k}_i} \cdots
~~~~\mbox{where}~~~ k_i ~=~\hat{k}_i + \tilde{k}_i.
\end{equation}
A factorization structure ~$\mathfrak{N}$ ~is {\bf embedded}\, in
~$\mathfrak{M}$, ~denoted by ~$\mathfrak{N} \prec \mathfrak{M}$, ~if
~$\mathfrak{N}=\mathfrak{M}$ ~or
~$\mathfrak{M}$ ~can be obtained by applying a sequence of embedding
operations on ~$\mathfrak{N}$.
~For example,
\begin{eqnarray*}
(4,3)^5 (1,6)^2 (3,2) & \prec & (4,3)^3 (4,3)^2 (1,6)^2 (3,2)
~~~~\mbox{(splitting ~$(4,3)^5$ ~to ~$(4,3)^3(4,3)^2$)} \\
& \prec & (4,3)^3 (5,9)^2 (3,2) ~~~~~~~~~~~~
\mbox{(combining ~$(4,3)^2 (1,6)^2$ ~to ~$(5,9)^2$)}.
\end{eqnarray*}
The relation ~$\prec$ ~is a partial ordering among factorization structures.

\vspace{-3mm}
\begin{theorem}[Factorization Manifold Embedding Theorem] \label{t:emb}
~Let
~$\mathfrak{N}$ ~be a factorization structure and
~$f \in \mathscr{F}^\mathfrak{N}$.
~For any factorization structure ~$\mathfrak{M}$ ~with
~$\dg(\mathfrak{M}) = \dg(\mathfrak{N})$ $= \mathbf{m}$, ~we have ~$f \in
\overline{\mathscr{F}^\mathfrak{M}}$ ~if and only if ~$\mathfrak{N}
\prec \mathfrak{M}$.
~Furthermore,
~$\codim(\mathscr{F}^\mathfrak{N}) \,>\, \codim(\mathscr{F}^\mathfrak{M})$
~in ~$\mathbb{P}^\mathbf{m}$
~if ~$\mathfrak{N} \prec \mathfrak{M}$ ~and ~$\mathfrak{N} \ne \mathfrak{M}$.
\end{theorem}

\vspace{-3mm}
\proof
~Assume ~$\mathfrak{N} \prec \mathfrak{M}$.
~To prove ~$ f \in \overline{\mathscr{F}^\mathfrak{M}}$,
~it suffices to show ~$f \in \overline{\mathscr{F}^{\tilde{\mathfrak{N}}}}$
~if ~$\tilde{\mathfrak{N}}$ ~is obtained from ~$\mathfrak{N}$
~by either one of the two embedding operations (\ref{dgcmb}) and
(\ref{msplt}).
~If ~$\tilde{\mathfrak{N}}$ ~is obtained by degree combining
(\ref{dgcmb}), then we can write ~$f = f_i^k f_j^k g$ ~with
~$f_i \in \mathscr{P}^{\mathbf{n}_i}$ ~and
~$f_j \in \mathscr{P}^{\mathbf{n}_j}$ ~being irreducible and coprime.
~By Corollary~\ref{c:od}, there is a polynomial sequence
~$\{h_l\}_{l=1}^\infty \subset
\mathscr{P}^{\mathbf{n}_i+\mathbf{n}_j}$ ~converging to zero such
that ~$f_i f_j + h_l$ ~is irreducible for all ~$l = 1, 2, \ldots$.
~Thus ~$f \in \overline{\mathscr{F}^{\tilde{\mathfrak{N}}}}$ ~since
~$(f_i f_j + h_l)^k g \in \mathscr{F}^{\tilde{\mathfrak{N}}}$
~converges to ~$f$ ~for
~$l \rightarrow \infty$.
~If ~$\tilde{\mathfrak{N}}$ ~is obtained by multiplicity splitting
(\ref{msplt}), then we can write ~$f = f_i^{k_i} g$ ~with
~$f_i \in \mathscr{P}^{\mathbf{n}_i}$.
~There is a sequence ~$\{h_l\}_{l=1}^\infty \subset
\mathscr{P}^{\mathbf{n}_i}$ ~converging to zero such that
~$f_i + h_l$ ~is irreducible for all ~$l = 1, 2, \ldots$.
~Thus ~$f \in \overline{\mathscr{F}^{\tilde{\mathfrak{N}}}}$ ~since
~$f_i^{\hat{k}_i}(f_i + h_l)^{\tilde{k}_i} g \in
\mathscr{F}^{\tilde{\mathfrak{N}}}$
~with ~$\hat{k}_i+\tilde{k}_i = k_i$
~converges to ~$f$ ~for ~$l \rightarrow \infty$.
~Conversely, assume ~$f \in \overline{\mathscr{F}^\mathfrak{M}}$
~with a irreducible squarefree  factorization
~$\alpha f_1^{k_1}\cdots f_r^{k_r}$.
~There is a sequence ~$\{g_l\}_{l=1}^\infty \subset \mathscr{F}^\mathfrak{M}$
~converging to ~$f$.
~Write ~$\mathfrak{M} = \mathbf{m}_1^{k_1^\prime}\cdots
\mathbf{m}_s^{k_s^\prime}$ ~and
~$g_l = \beta_l g_{l1}^{k_1^\prime}\cdots g_{ls}^{k_s^\prime}$ ~for
~$l=1,2,\ldots$.
~By Lemma~\ref{l:pcv}, we can further assume
~$\displaystyle \lim_{l\rightarrow\infty} g_{lj} =
\hat{g}_j \in \mathscr{P}^{\mathbf{m}_j}$ ~for ~$j=1,\ldots,s$
~and ~$\beta_l \rightarrow \hat\beta$.
~Due to ~$\alpha f_1^{k_1}\cdots f_r^{k_r} = \hat\beta
\hat{g}_{1}^{k_1^\prime}\cdots \hat{g}_{s}^{k_s^\prime}$ ~and the
uniqueness of factorizations, we can factor polynomials
~$\hat{g}_1,\ldots,\hat{g}_s$ ~and
combine equivalent irreducible factors into higher multiplicities
to reproduce
the squarefree irreducible factorization ~$\alpha f_1^{k_1}\cdots f_r^{k_r}$.
~Namely, the structure ~$\mathfrak{M}$ ~can be obtained by a sequence of
embedding operations on ~$\mathfrak{N}$, ~leading to
~$\mathfrak{N}\prec\mathfrak{M}$.
~The inequality
~$\codim(\mathscr{F}^\mathfrak{N}) \,>\, \codim(\mathscr{F}^\mathfrak{M})$
~follows from a straightforward verification
using (\ref{eq:codim}) on (\ref{dgcmb}) and (\ref{msplt}).
\foorp

\begin{figure}
\begin{center}
\epsfig{figure=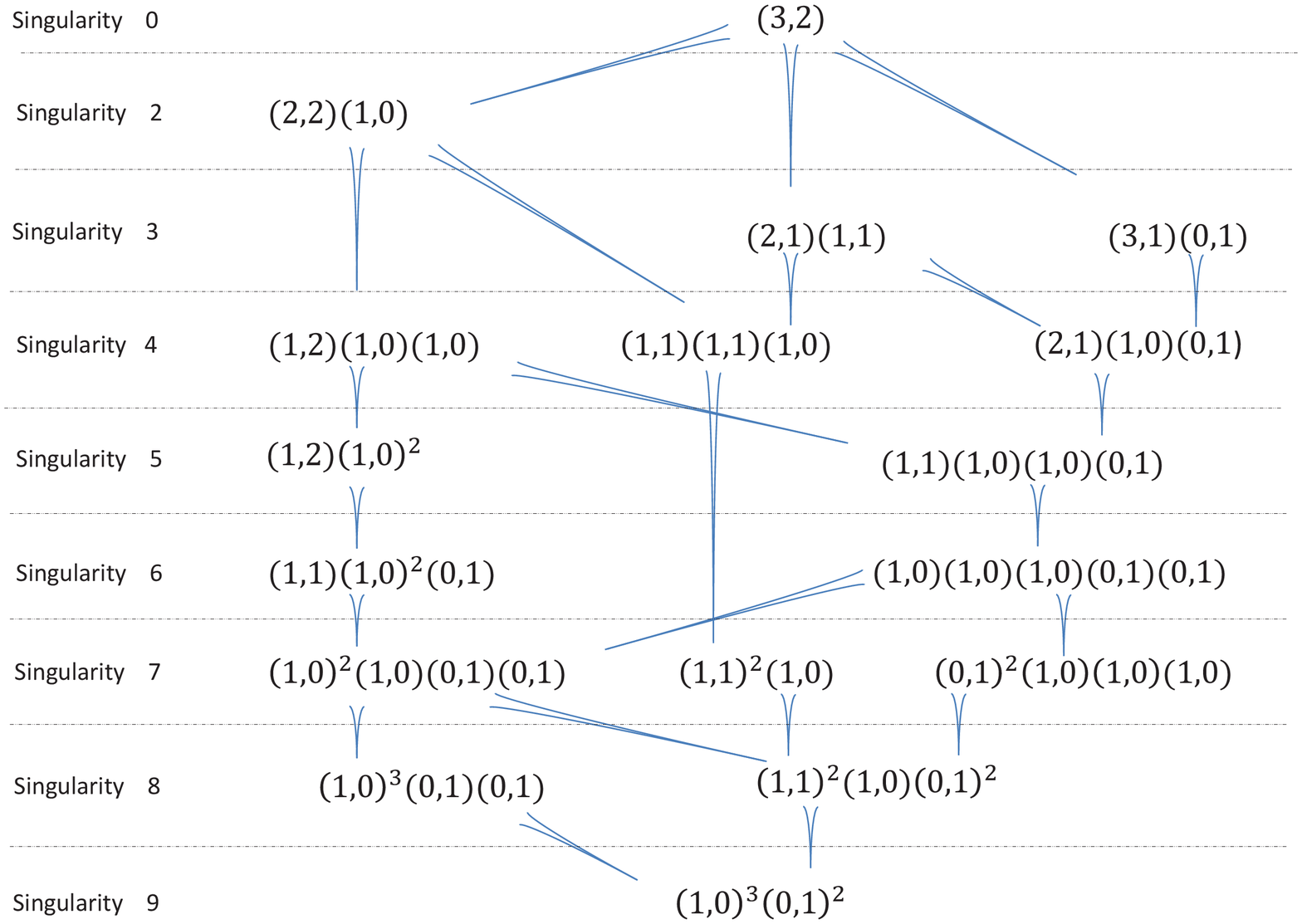,width=4.5in,height=2.8in}
   \vspace{-5mm}
   \caption{\small Stratification of factorization manifolds in
~$\mathscr{P}^{(3,2)}$, ~where ~$\mathfrak{N} \prec \mathfrak{M}$
~indicates ~$\mathscr{F}^\mathfrak{N} \subset
\overline{\mathscr{F}^\mathfrak{M}}$.}
   \label{fig:1}
\end{center}
\vspace{-6mm}
\end{figure}

The Factorization Manifold Embedding Theorem
implies the geometry of polynomial factorization:
~{\em The subset
~$\mathscr{P}^\mathbf{m}$ ~of degree ~$\mathbf{m}$ ~polynomials is a
disjoint union of factorization manifolds that are
topologically stratified in such a way that every factorization manifold
of positive singularity is embedded in the closure of a
factorization manifold of lower singularity.}
~As an example, Figure \ref{fig:1} illustrates such a stratification among
all the factorization manifolds through corresponding factorization
structures in ~$\mathscr{P}^{(3,2)}$.

We define the distance between a polynomial and a factorization
manifold
\begin{equation} \label{dfF}
\dist\big(f, \mathscr{F}^{\mathfrak{M}}\big) ~~=~~ \inf_{g\in
\mathscr{F}^{\mathfrak{M}}} (||f - g||).
\end{equation}
Let ~$\mathfrak{N}$ ~be the factorization structure of ~$f$.
~The distance ~$\dist\big(f, \mathscr{F}^{\mathfrak{M}}\big) = 0$
~if and only if
~$f \in \overline{\mathscr{F}^{\mathfrak{M}}}$, ~which
is equivalent to ~$\mathfrak{N} \prec \mathfrak{M}$ ~by the
Factorization Manifold Embedding Theorem.
~As a consequence, the native manifold ~$\mathscr{F}^\mathfrak{N}$ ~of ~$f$
~distinguishes itself as the unique factorization manifold that is of the
{\em highest singularity} (i.e. highest codimension) among all the
factorization manifolds having a distance zero to ~$f$.
~More precisely, a polynomial ~$f$ ~belongs to a factorization manifold
~$\mathscr{F}^{\mathfrak{N}}$ ~if and only if, 
in ~$\mathbb{P}^\mathbf{m}$,
\[ \codim (\mathscr{F}^{\mathfrak{N}}) ~~=~~ \max\Big\{\codim
(\mathscr{F}^{\mathfrak{M}}) ~\Big|
~\dg(\mathfrak{M}) = \dg(\mathfrak{N}) =\mathbf{m} ~~\mbox{and}~~
\dist\big(f, \mathscr{F}^{\mathfrak{M}}\big)=0   \Big\}. \]

On the other hand, a polynomial ~$\tilde{f} \in
\mathscr{F}^{\mathfrak{N}}$ ~with
~$\mathfrak{N} \not\prec \mathfrak{M}$ ~implies
~$\dist\big(\tilde{f}, \mathscr{F}^{\mathfrak{M}}\big) > 0$.
~Since there are finitely many factorization manifolds, there
exists a minimum positive distance
\begin{equation} \label{theta}
\theta_{\tilde{f}} ~~=~~ 
\min_{\overline{\mathscr{F}^\mathfrak{M}} \,\not\ni\, \tilde{f}}
\,\dist\big(\tilde{f}, \mathscr{F}^{\mathfrak{M}}\big) ~~>~~ 0.
\end{equation}
The constant ~$\theta_{\tilde{f}}$ ~is the {\em critical gap} of 
~$\tilde{f}$ ~from unembedded singularities and it is the very
window of opportunity for numerical factorization.
~When the polynomial ~$\tilde{f} \in
\mathscr{F}^\mathfrak{N}$ ~is represented by
an empirical version ~$f$ ~with a small perturbation
~$\|f-\tilde{f}\| < \frac{1}{2}\,\theta_{\tilde{f}}$,
~the underlying factorization
structure can still be identified by the following lemma.

\vspace{-3mm}
\begin{lemma}\label{lem:maxcodim}
~Let ~$\tilde{f}$ ~be a polynomial with a factorization structure
~$\mathfrak{N}$ ~with ~$\theta_{\tilde{f}}$ ~be given in {\em (\ref{theta})}.
~For any empirical data ~$f$ of ~$\tilde{f}$ ~satisfying
~$\|f-\tilde{f}\| < \frac{1}{2}\, \theta_{\tilde{f}}$, ~the factorization
structure ~$\mathfrak{N}$ ~of ~$\tilde{f}$ ~is uniquely identifiable
using the data ~$f$ ~by
\begin{equation} \label{mcde}
 \codim (\mathscr{F}^{\mathfrak{N}}) ~~=~~ \max \big\{\codim
(\mathscr{F}^{\mathfrak{M}}) \mid
~\dg(\mathfrak{M}) = \dg(\mathfrak{N}) =\mathbf{m}~~\mbox{and}~~
\dist\big(f, \mathscr{F}^{\mathfrak{M}}\big) < \epsilon   \big\}
\end{equation}
in ~$\mathbb{P}^\mathbf{m}$ ~for any ~$\epsilon$ ~satisfying
~$\dist\big(f, \mathscr{F}^\mathfrak{N}\big) < \epsilon <
\frac{1}{2}\,\theta_{\tilde{f}}$.
\end{lemma}

\vspace{-3mm}
\proof
A straightforward verification.
\foorp

In summary, singular polynomials form factorization
manifolds with positive codimensions and nonsingular polynomials
form an open dense subset in ~$\mathbb{P}^\mathbf{m}$.
~Those factorization manifolds topologically stratify in such a way that
every singular manifold belongs to the
closures of some manifolds of lower singularities.
~Almost all tiny perturbations on a singular polynomial alter
its factorization structure in such a way that the singularity
reduces and never increases.
~There is a gap from any singular polynomial to higher singularity and
this gap ensures the lost factorization structure can be recovered by
finding the highest singularity manifold nearby if
the perturbation is small.
~As a result, identifying the factorization structure is  well-posed
as an optimization problem.

\vspace{-5mm}
\section{The notion of numerical factorization} \label{s:nf}

\vspace{-5mm}
We shall rigorously formulate the concept of the numerical factorization to
remove the ill-posedness of the conventional factorization, and to achieve
the main objective of recovering the exact factorization accurately using
the imperfect empirical data.
~The numerical factorization should approximate the underlying factorization
with an accuracy the data deserve.
~The following problem statement gives a precise description of the
problem that numerical factorization is intended to solve.

\vspace{-3mm}
\begin{Problem}[Numerical Factorization Problem]\label{pro:nif}
~Let ~$f$ ~be a polynomial as the empirical data of
an underlying polynomial ~$\tilde{f}$ ~whose
irreducible factorization ~$\tilde{\alpha}
\tilde{f}_1 \cdots \tilde{f}_r$ ~is to be computed.
~Assuming the data error ~$\big\|f - \tilde{f}\big\|$ ~is sufficiently
small, find an irreducible
factorization ~$\alpha f_1 \cdots f_r$
~of a certain polynomial ~$\hat{f}$ ~such that both the backward error
and forward error are in the order of data error and the unit round-off:
\begin{eqnarray} 
\big\|f - \alpha f_1 \cdots f_r\big\| & = & 
O \big(\|f- \tilde{f}\|+\mathfrak{u}\big) \label{bdoff} \\
\dist\big(\alpha f_1 \cdots f_r, 
~\tilde{\alpha} \tilde{f}_1 \cdots \tilde{f}_r
\big)
 &~~=~~& O \big(\|f- \tilde{f}\|+\mathfrak{u}\big).
\label{doff}
\end{eqnarray}
where ~$\mathfrak{u}$ ~is the unit round-off in the floating point arithmetic.
\end{Problem}

\vspace{-3mm}
Notice that (\ref{doff}) implies 
~$\alpha f_1 \cdots f_r$ ~and 
~$\tilde{\alpha} \tilde{f}_1 \cdots \tilde{f}_r$
~are required to have the same factorization structure
by the definition of the distance (\ref{ferr}).
~Problem~\ref{pro:nif} goes a step further from the Open Problem 1 in 
\cite{kal99} in which only the backward error is required to be small.

\begin{figure}[ht]
\begin{center}
   \epsfig{figure=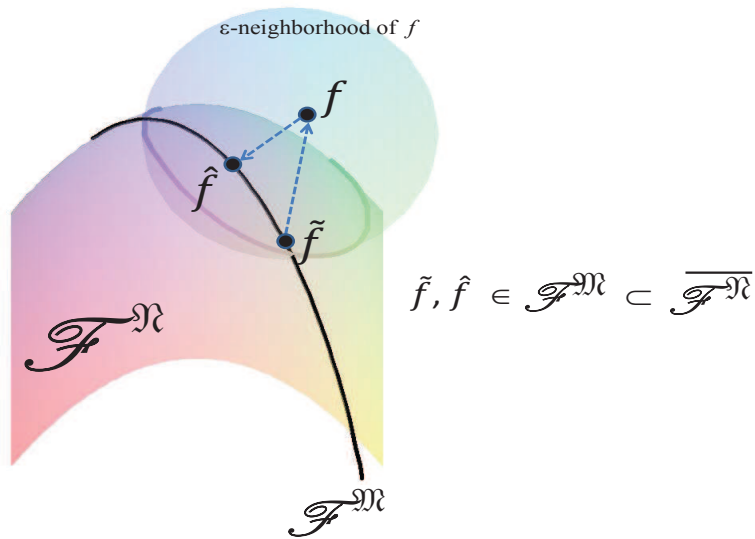,width=5.0in,height=3.2in}
    \vspace{-20pt}
   \caption{\small Illustration of the numerical factorization of
~$f$}
   \label{fig:2}
\end{center}
\vspace{-8mm}
\end{figure}

Let ~$\tilde{f}$ ~be the polynomial in Problem~\ref{pro:nif}
with the factorization structure ~$\mathfrak{M}$
~and ~$f$ ~be its empirical data
representation, as illustrated in Figure~\ref{fig:2}.
~By the Factorization Manifold Embedding Theorem, the data polynomial ~$f$ 
~is away
from the native factorization manifold ~$\mathscr{F}^\mathfrak{M}$ ~with a
{\em reduced singularity}.
~Note that the data ~$f$ ~is also near all the factorization
manifolds ~$\mathscr{F}^\mathfrak{N}$ ~with ~$\mathfrak{M} \prec
\mathfrak{N}$ ~and the native manifold ~$\mathscr{F}^\mathfrak{M}$
~is not the nearest in distance but highest in singularity by
Lemma~\ref{lem:maxcodim}.
~Upon identifying the factorization structure ~$\mathfrak{M}$, ~it
is then natural to calculate the exact irreducible factorization
~$\alpha f_1^{k_1} \cdots f_r^{k_r}$ ~of the polynomial ~$\hat{f}
\in \mathscr{F}^\mathfrak{M}$ ~that is the nearest to ~$f$ ~and
designate it as the numerical irreducible factorization of
~$f$ ~since Corollary~\ref{c:iflip} suggests that the (exact)
irreducible factorization of ~$\hat{f}$ ~approximates that of
~$\tilde{f}$.
~The following definition is the detailed formulation.

\vspace{-3mm}
\begin{define}[Numerical Factorization]\label{def:NIF}
~For a given polynomial ~$f$ ~and a backward error tolerance
~$\epsilon > 0$, we say ~$\alpha f_1 f_2 \cdots f_s$ ~is a
{\bf numerical irreducible factorization of ~$f$ ~within
~$\epsilon$} ~if
~$\alpha f_1 f_2 \cdots f_s \in \mathscr{F}^\mathfrak{M}$ ~is
an irreducible factorization and
\begin{equation} \label{mdis}
\min_{\gamma \in \C}\big\|f - \gamma f_1 f_2 \cdots f_s\big\|
~~=~~ \min_{g \in\mathscr{F}^\mathfrak{M}} \|f-g\|
~~=~~ \dist\big(f, \,\mathscr{F}^\mathfrak{M}\big) ~~<~~ \epsilon,
\end{equation}
where ~$\mathfrak{M}$ ~is the factorization structure of the degree
~$\mathbf{m} = \dg(f)$ ~such that
\begin{equation} \label{mcod}
\codim(\mathscr{F}^\mathfrak{M}) ~~=~~
\max \big\{
\codim (\mathscr{F}^\mathfrak{N} ) ~\big|
~\dg(\mathfrak{N}) = \mathbf{m} ~~\mbox{and}~~
\dist(f,\mathscr{F}^\mathfrak{N}) < \epsilon \big\}
\end{equation}
in ~$\mathbb{P}^{\mathbf{m}}$.
~We call ~$\alpha f_1^{k_1}\cdots f_r^{k_r}$ ~a
{\bf numerical irreducible squarefree factorization of ~$f$ ~within
~$\epsilon$} ~if it is squarefree and it is a numerical
irreducible factorization of ~$f$ ~within ~$\epsilon$.
\end{define}

\vspace{-3mm}
We shall use the abbreviated term {\bf numerical factorization} for either
the
numerical irreducible factorization or the
numerical irreducible squarefree
factorization when the distinction is insignificant in the context.
~The formulation of the numerical factorization follows the same
``three-strikes'' principles that have been effectively applied to
the regularization of other ill-posed algebraic problems
\cite{Zeng09b}:
~The numerical factorization
of ~$f$ ~is the exact factorization of a nearby polynomial
~$\hat{f}$ ~within a backward error tolerance ~$\epsilon$ (backward
nearness principle).
~The nearby polynomial ~$\hat{f}$ ~is of the highest singularity among all
the polynomials in the ~$\epsilon$-neighborhood of ~$f$ (maximum singularity
principle).
~The nearby polynomial ~$\hat{f}$ ~is the
nearest polynomial to the given ~$f$ ~among all the polynomials
with the same singularity as ~$\hat{f}$ (minimum distance principle).

The error tolerance ~$\epsilon$ ~in Definition \ref{def:NIF} ~depends on the
particular application, the hardware precision, the underlying
polynomial ~$\tilde{f}$ ~and the data error ~$\|f-\tilde{f}\|$.
~The interval for setting ~$\epsilon$ ~will be established in the
Numerical Factorization Theorem in \S\ref{s:reg}.
~Notice that a polynomial ~$f$ ~can easily have different numerical
factorizations within different error tolerances approximating
different factorizations
(c.f. Example~\ref{ex:5f} in \S\ref{s:imp}).
%

\vspace{-5mm}
\section{Regularity and sensitivity of numerical factorization}
\label{s:reg}

\vspace{-5mm}
As a concept attributed to Jacques S. Hadamard,
a mathematical problem is well-posed if its
solution holds existence, uniqueness and continuity with respect 
to data.
~Furthermore, Lipschitz continuity of the solution is crucial for 
numerical computation as it implies a finite sensitivity with respect to
data perturbations and round-off.
%
%
~With the geometry established in \S\ref{s:geo} and \S\ref{s:str} ,
the well-posedness of numerical factorization is a direct
consequence of the Tubular Neighborhood Theorem, which is
one of the fundamental
results in differential topology.
~The following elementary version of the Tubular Neighborhood Theorem is
adapted from its abstract form for complex analytic manifolds in ~$\C^n$.

\vspace{-3mm}
\begin{lemma}[Tubular Neighborhood Theorem]\label{l:tnt}
{\em \cite{ZZ-tnt}}
~Every complex analytic manifold is contained in a tubular neighborhood.
~More precisely,
~for every complex analytic manifold ~$\Pi$ ~in ~$\C^n$, ~there is an open subset
~$\Omega$ ~of ~$\C^n$ ~containing ~$\Pi$ ~and a projection mapping
~$\pi ~:~ \Omega \longrightarrow \Pi$ ~such that, for every ~$\mathbf{z}
\in \Omega$, ~its projection ~$\pi(\mathbf{z}) \in \Pi$ ~is the
unique distance-minimization point from ~$\mathbf{z}$ ~to ~$\Pi$, ~namely
~$\|\pi(\mathbf{z}) - \mathbf{z}\|_2 ~=~ \displaystyle
\min_{\mathbf{u} \in \Pi} \big\|\mathbf{u} - \mathbf{z}\big\|_2$.
~Furthermore, the mapping ~$\pi$ ~is locally Lipschitz continuous.
\end{lemma}

\vspace{-2mm}
We can now establish the main theorem, which 
asserts the properties that are desirable from the numerical
factorization as formulated in Definition~\ref{def:NIF}, provides a complete 
regularization and, in essence, achieves the objectives of numerical
factorization in Problem~\ref{pro:nif}.

\vspace{-3mm}
\begin{theorem}[Numerical Factorization Theorem]
\label{th:regular}
~Let ~$\tilde{f}$ ~be a polynomial of degree ~$\mathbf{m}$ ~with its critical 
gap ~$\theta_{\tilde{f}}$ ~as in {\em (\ref{theta})} and
an irreducible factorization 
~$\tilde\alpha\tilde{f}_1 \tilde{f}_2 \cdots \tilde{f}_k$.
~Then ~$\theta_{\tilde{f}} > 0$ ~and the following properties of numerical
factorization hold.
\vspace{-3mm}
\begin{itemize}\parskip-2mm
\item[\em (i)] ~{\em Conventional factorization is a special case of 
numerical factorization:} ~The numerical factorization of ~$\tilde{f}$ ~within
any ~$\epsilon\in (0,\theta_{\tilde{f}})$ ~is identical to the exact irreducible 
factorization of ~$\tilde{f}$.
\item[\em (ii)] ~{\em Computing numerical factorization is a well-posed 
problem:} ~There is a neighborhood ~$\Omega_{\tilde{f}}$ ~of
~$\tilde{f}$ ~in ~$\mathbb{P}^\mathbf{m}$
~such that every ~$f \in \Omega_{\tilde{f}}$ ~is associated with a
constant ~$\delta_f \le \|f-\tilde{f}\|$ ~such that the numerical factorization
~$\alpha f_1 f_2 \cdots f_l$ ~of ~$f$ ~uniquely exists within ~$\epsilon$ 
~for all ~$\epsilon \in \left(\delta_f, \frac{1}{2}\,\theta_{\tilde{f}}\right)$
~and is Lipschitz continuous with respect to ~$f$.
\item[\em (iii)] ~{\em Numerical factorization is backward accurate:} 
~For every ~$f\in\Omega_{\tilde{f}}$ ~and ~$\epsilon \in \left(\delta_f, 
\frac{1}{2}\,\theta_{\tilde{f}}\right)$, ~the numerical factorization
~$\alpha f_1 f_2 \cdots f_l$ ~of ~$f$ ~within ~$\epsilon$ ~satisfies
\begin{equation} 
\big\|f-\alpha f_1 f_2 \cdots f_l\big\| ~~\le~~  \|f-\tilde{f}\|.
\label{bdftf} \end{equation}
\item[\em (iv)] ~{\em The conventional factorization can be accurately 
recovered from empirical data:}
~For every ~$f\in\Omega_{\tilde{f}}$ ~as empirical data of ~$\tilde{f}$
and ~$\epsilon \in \left(\delta_f, 
\frac{1}{2}\,\theta_{\tilde{f}}\right)$, ~the numerical factorization
~$\alpha f_1 f_2 \cdots f_l$ ~of ~$f$ ~within ~$\epsilon$ 
~has the identical structure as 
~$\tilde\alpha \tilde{f}_1 \tilde{f}_2 \cdots \tilde{f}_k$ 
~and 
\begin{equation} 
 \dist\big(\alpha f_1 f_2 \cdots f_l,
~\tilde\alpha \tilde{f}_1 \tilde{f}_2 \cdots \tilde{f}_k
\big)  ~~<~~  \eta_{\tilde{f}} \, \|f-\tilde{f}\|.
\label{dftf}
\end{equation}
where ~$\eta_{\tilde{f}} >0$ ~is a constant depends on ~$f$. 
\end{itemize}
\end{theorem}

\vspace{-3mm}
\proof
Let ~$\mathfrak{M}$ ~denote the factorization structure of ~$\tilde{f}$.
~Then ~$\mathscr{F}^\mathfrak{M}$ ~is the manifold of the highest singularity 
within ~$\epsilon\,\in\,(0,\theta_{\tilde{f}})$ ~of ~$\tilde{f}$, ~and
~$\tilde{f}$ ~itself is the polynomial of minimum distance zero on 
~$\mathscr{F}^\mathfrak{M}$ ~from ~$\tilde{f}$, ~and thus (i) holds.
%
%
~Let ~$\Sigma$ ~be the tubular neighborhood of
~$\mathscr{F}^\mathfrak{M}$ ~described in Lemma~\ref{l:tnt} and 
let ~$\Omega_{\tilde{f}} \subset \Sigma$ ~be a neighborhood of ~$\tilde{f}$
~such that every ~$f \in \Omega_{\tilde{f}}$ ~satisfies
~$\|f-\tilde{f}\| <\frac{1}{2}\,\theta_{\tilde{f}}$.
~Set ~$\delta_f = \dist\big(f,\mathscr{F}^\mathfrak{M}\big)$.
~Then, for every ~$\epsilon \in (\delta_f,\,\frac{1}{2}\theta_{\tilde{f}})$
~the equality (\ref{mcod}) holds
since ~$\dist\big(f,\mathscr{F}^\mathfrak{N}\big) > \frac{1}{2}\,
\theta_{\tilde{f}} > \epsilon$ ~for every ~$\mathfrak{N} \not\succ \mathfrak{M}$.
~By the Tubular Neighborhood Theorem, there exists a unique
~$\hat{f} = \pi(f) \in \mathscr{F}^\mathfrak{M}$ ~with minimal distance to
~$f$.
~As a result, the numerical factorization of ~$f$ ~uniquely exists
as the exact irreducible factorization of ~$\hat{f}$, ~and the
numerical factorization is locally Lipschitz continuous since
~$\pi$ ~is locally Lipschitz continuous along with
Corollary~\ref{c:iflip}, leading to part (ii).
~Part (iii) is true since ~$\|f-\hat{f}\| \le \|f-\tilde{f}\|$.
~The Lipschitz continuity of the numerical factorization
also implies (\ref{dftf}) and part (iv). \foorp

In simpler terms, Numerical Factorization Theorem ensures that every
factorable polynomial ~$\tilde{f}$ ~is allowed to be perturbed while
its factorization can still be recovered as long as the empirical
data ~$f$ ~is still in the neighborhood
~$\Omega_{\tilde{f}}$.
~For each data representation ~$f$ ~of ~$\tilde{f}$, ~there is a window
~$(\delta_f, \frac{1}{2}\,\theta_{\tilde{f}})$ ~for setting the error tolerance
~$\epsilon$ ~for recovering the factorization of ~$\tilde{f}$.
~The fact that the lower bound ~$\delta_f$ ~of the error tolerance 
~$\epsilon$ ~is no larger than the data error ~$\|f-\tilde{f}\|$ ~is 
significant in practical computation:
~If a data error bound ~$\eta > 0$ ~for ~$\|f-\tilde{f}\|$ ~is known
or can be estimated in an application, the error tolerance can be set at 
~$\epsilon = \eta$ ~or a moderate multiple of the unit round-off, 
whichever is larger.
~The upper bound ~$\frac{1}{2}\,\theta_{\tilde{f}}$ ~appears to be difficult
to estimate but not needed as long as it is not too small.

With a proper error tolerance ~$\epsilon$, ~the numerical factorization of 
the data ~$f$ ~within ~$\epsilon$
~approximates the exact factorization of the underlying polynomial
~$\tilde{f}$ ~with an accuracy in the same order of the data accuracy.
~Namely, the numerical factorization we formulated in Definition~\ref{def:NIF}
achieves the objective of the numerical factorization problem as specified
in Problem~\ref{pro:nif}.
~Furthermore, computing the numerical factorization is a well-posed
problem with a finite sensitivity that can be established in the following
theorem.

\vspace{-3mm}
\begin{theorem}[Numerical Factorization Sensitivity Theorem]
~Let ~$\alpha f_1^{k_1} \cdots f_r^{k_r}$ ~be the numerical factorization
of ~$f$ ~within certain ~$\epsilon$ ~and ~$g$ ~be sufficiently close to ~$f$
~so that its numerical factorization within ~$\epsilon$ ~can be written as
and ~$\gamma g_1^{k_1} \cdots g_r^{k_r}$ ~with
~$\dg(g_j) = \dg(f_j) = \mathbf{m}_j$ ~for ~$j=1,2,\ldots,r$.
~Further assume
~$\Jac(\cdot)$ ~is as defined in {\em (\ref{eq:Jac_diffeo2})} where
~$\mathbf{b}_j \in \C^{\langle\mathbf{m}_j\rangle}$ ~with
~$\|\mathbf{b}_j\|_2 = 1$ ~and
~$\mathbf{b}_j \cdot \lb f_j \rb =1$ ~for ~$j=1,\ldots,r$.
~Then
\begin{equation} \label{nifss}
\limsup_{g \rightarrow f}
\frac{\dist\big(\alpha f_1^{k_1} \cdots f_r^{k_r},
~\gamma g_1^{k_1} \cdots g_r^{k_r}\big)
}{\|f-g\|} ~~\le~~ \eta\,
\Big\|\Jac(\alpha,\lb f_1 \rb, \ldots, \lb f_r \rb)^+\Big\|_2
~~<~~ \infty
\end{equation}
where ~$\eta$ ~is a constant associated with ~$f$ ~and ~$\mathbf{b}_1,
\ldots,\mathbf{b}_r$.
\end{theorem}

\vspace{-3mm}
\proof
~A straightforward verification using Theorem 2 in
\cite{ZZ-tnt} and Lemma~\ref{l:dis}.
\foorp

The inequality (\ref{nifss}) depends on the choices of the specific
representative ~$\alpha f_1^{k_1} \cdots f_r^{k_r}$ ~in the
equivalent class of factorizations and the scaling vectors
~$\mathbf{b}_1,\ldots,\mathbf{b}_r$.
~Independent of those choices, we define the positive real number
\begin{eqnarray}
\lefteqn{\kappa_{\epsilon}(f) ~~:=~~ \inf \Big\{
\,\big\|\Jac(\beta,\lb h_1 \rb, \ldots, \lb h_r \rb)^+\big\|_2 ~\Big|~
\beta h_1^{k_1}\cdots h_r^{k_r} \sim \alpha f_1^{k_1} \cdots f_r^{k_r},}
\nonumber \\\label{nfcond}
& & ~~~~~~~~~~~~~~~~~~~~~~~~~~~
\mathbf{b}_j \in \C^{\langle\mathbf{m}_j\rangle},
~\|\mathbf{b}_j\|_2 = 1,
~\mathbf{b}_j \cdot \lb h_j \rb =1, ~j=1,\ldots,r
\Big\}
\end{eqnarray}
as the {\bf condition number of the numerical factorization}
~of ~$f$ ~within ~$\epsilon$ ~where
~$\alpha f_1^{k_1} \cdots f_r^{k_r}$ ~is a numerical squarefree irreducible
factorization of ~$f$ ~within ~$\epsilon$.
~From Lemma~\ref{lem:injective}, this condition number is finite
since the factorization ~$\alpha f_1^{k_1} \cdots f_r^{k_r}$ ~of ~$f$ ~is
squarefree, and ~$\kappa_{\epsilon}(f)$ ~becomes large when the Jacobian
~$\Jac(\beta,\lb h_1 \rb, \ldots, \lb h_r \rb)^+$  ~is near rank-deficient
when two of the factors ~$h_1,\ldots,h_r$ ~are a small perturbation away
from having nonconstant GCD.
~Consequently, the nature of the computation stability of numerical
factorization become apparent:
{\em
~The numerical factorization
~$\alpha f_1^{k_1} \cdots f_r^{k_r}$ ~of ~$f$
~is ill-conditioned if there exist
two factors ~$f_i$ ~and ~$f_j$ ~that are near
non-coprime polynomials so
that a small perturbation of ~$f$ ~can increase
the singularity above that of ~$\hat{f} = \alpha f_1^{k_1} \cdots f_r^{k_r}$.
}

\section{On the numerical squarefree factorization} \label{s:nsf}

\vspace{-5mm}
Every polynomial ~$f$ ~has a unique squarefree factorization
~$\alpha f_1^{k_1}\cdots f_r^{k_r}$ ~where
~$f_1,\ldots,f_k$ ~are pairwise coprime squarefree polynomials.
~Such squarefree factorizations are important in its own right and usually
easier to compute than irreducible factorizations.
~Our numerical factorization algorithm and implementation start with
finding a numerical squarefree factorization followed by numerical
irreducible factorizations of the squarefree components.
~Naturally, the notion of numerical squarefree factorization and
its properties are in question.

Similar to (irreducible) factorization structure, we can define a
{\em squarefree factorization structure}
~$\mathfrak{N} = \mathbf{n}_1^{k_1}\cdots\mathbf{n}_r^{k_r}$ ~of polynomials
having a squarefree factorization
~$f = \alpha f_1^{k_1}\cdots f_r^{k_r}$ ~where ~$k_1\le\cdots\le k_r$ ~where
~$f_j\in\mathbb{P}^{\mathbf{n}_j}$ ~is squarefree for ~$j=1,\ldots,r$ ~and
pairwise coprime.
~Also let ~$\mathscr{S}^\mathfrak{N}$ ~denote the collection of polynomials
in ~$\mathbb{P}^\mathbf{m}$ ~having a squarefree factorization structure
~$\mathfrak{N}$.
~Notice that Lemma~\ref{lem:injective} applies to squarefree factorizations
since the irreducibility of factors is not required.
~It is also a straightforward verification that Lemma~\ref{l:pcv} and
Corollary~\ref{c:phi} still hold for ~$\mathscr{S}^\mathfrak{N}$.
~As a result, the subset ~$\mathscr{S}^\mathfrak{N}$ ~is also a complex analytic
manifold in ~$\mathbb{P}^\mathbf{m}$ ~of codimension
\[
 \codim (\mathscr{S}^\mathfrak{N}) ~~=~~
\langle \mathbf{m} \rangle
 - \Big(\langle \textbf{n}_{1} \rangle + \cdots +
\langle \textbf{n}_{r} \rangle +1 - r\Big).
\]
The embedding properties of squarefree factorization manifolds hold as well.

Similar to Definition~\ref{def:NIF}, we can formulate the
{\em numerical squarefree factorization} of a polynomial
~$f$ ~within an error tolerance ~$\epsilon$ ~as the exact factorization of a
polynomial ~$\hat{f}\in\mathscr{S}^\mathfrak{N}$ ~where
~$\mathscr{S}^\mathfrak{N}$ ~is the squarefree factorization manifold of the
highest codimension among all manifolds intersecting the
~$\epsilon$-neighborhood of ~$f$ ~and ~$\hat{f}$ ~is the nearest polynomial
from ~$f$ ~on ~$\mathscr{S}^\mathfrak{N}$.
~Such a numerical squarefree factorization is a generalization of the
conventional exact squarefree factorization and accurate approximation to
the exact squarefree factorization of the underlying polynomial ~$\tilde{f}$.
~Furthermore its computation is a well-posed problem with a finite
sensitivity measure.

For a unit vector ~$\mathbf{z} = (z_1,\ldots,z_{\ell}) \in
\C^\ell$, ~let ~$\partial_\mathbf{z}p$ ~denote the directional
derivative of ~$p \in \C[x_1,\ldots,x_\ell]$ ~along (the direction
of) ~$\mathbf{z}$, ~namely ~$\partial_\mathbf{z} p \,=\,
z_1\, \frac{\partial p}{\partial x_1} + \cdots
+z_\ell\, \frac{\partial p}{\partial x_\ell}$.
~The following lemma is the basis for the numerical
squarefree factorization.

\vspace{-3mm}
\begin{lemma}
~Every ~$f \in \C[x_1,\ldots,x_\ell]$ ~has a squarefree factorization
~$\alpha h_1^{k_1}\cdots h_r^{k_r}$ ~with non-constant factors
~$h_1,\ldots,h_r$ ~and distinct multiplicities ~$k_1,\ldots,k_r \ge 1$.
~Furthermore, for almost all unit vectors ~$\mathbf{z} \in \C^\ell$,
\begin{equation} \label{hj}
\GCD(v,\,l \cdot \partial_\mathbf{z}v-w) ~~=~~
\left\{ \begin{array}{ccl} h_j && ~\mbox{for}~~ l = k_j
~~~\mbox{with}~~ j = 1,\,\ldots\,,r \\
1 &&~\mbox{if}~~~ l \not\in \{k_1,\ldots,k_r\}
\end{array} \right.
\end{equation}
where ~$v$ ~and ~$w$ ~are cofactors of
~$u = \GCD(f,\partial_\mathbf{z} f)$
~such that ~$f=uv$ ~and ~$\partial_\mathbf{z}f = uw$.
\end{lemma}

\vspace{-3mm}
\proof
~The existence of ~$h_1,\ldots,h_r$ ~is obvious.
~For almost all ~$\mathbf{z} \in \C^\ell$, ~
{$\partial_\mathbf{z} h_j \,\ne\, 0$ ~for ~$j=1,\ldots,r$.
~Thus ~$v \,\sim\, h_1\cdots
h_r$ ~and ~$w ~\sim~
\displaystyle \sum_{j=1}^r k_j (\partial_\mathbf{z} h_j)\,
\prod_{i\ne j} h_i$},  ~leading to (\ref{hj}). \foorp

The numerical squarefree factorization can be computed by a sequence of
numerical greatest common divisor replacing the exact GCD in (\ref{hj}).

By Lemma~\ref{lem:injective}, the Jacobian of the mapping
~$\phi$ ~in
(\ref{eq:diffeo2}) is injective at the least squares solution of ~$
{\phi}(\beta,\lb g_1 \rb, \ldots, \lb g_r \rb) = (\lb f\rb,1,\ldots,1)$,
~implying the Gauss-Newton iteration locally converges to this
least squares solution if the data ~$f$ ~and the
initial iterate are sufficiently accurate.

Due to its similarity with the numerical irreducible squarefree factorization,
we omit the detailed elaboration of the numerical squarefree factorization
in this paper.

\vspace{-5mm}
\section{Computation of numerical factorizations}

\vspace{-5mm}
Overall, computing the numerical factorization consists of
two stages.
~The first stage identifies the factorization structure
along with initial approximations of the numerical factors.
~In the second stage,
the numerical factors are refined to minimize the distance from the given
polynomial to the manifold associated with the factorization structure.

In the first stage, the factorization structure can be computed by a sequence
numerical squarefree factorizations,
rank-revealing of the Ruppert matrices \cite{Gao03,kal-may,rup99},
generalized eigenvalue computation
and numerical greatest common divisor calculation.
~Initial approximations of the numerical factors are obtained
as by-products.

In the second stage, the initial factor approximations
~$p_1,\ldots,p_r$ ~of degrees ~$\mathbf{m}_1,\ldots,\mathbf{m}_r$ ~can be
scaled to unit norms ~$\big\| p_1 \big\| = \cdots = \big\|
p_r \big\| = 1$.
~The mapping ~$\phi$ ~in (\ref{eq:diffeo2}) becomes well defined by setting
up the scaling vectors ~$\mathbf{b}_i \,=\, \lb p_i \rb$ ~for ~$i = 1, \ldots, r$,
~The second stage of the numerical factorization algorithm is
essentially the process of solving for the least squares solution to
the overdetermined nonlinear system
\begin{equation} \label{nls}
\phi(\mathbf{z}) ~~=~~ (\lb f\rb, 1, \ldots, 1)
,
~~~~\mbox{with}~~~
\mathbf{z} \in \C\times \C^{\langle \mathbf{m}_1\rangle} \times \cdots \times
\C^{\langle \mathbf{m}_r\rangle}
\end{equation}
using the Gauss-Newton iteration
\begin{equation} \label{gni}
\mathbf{z}_{j+1} ~~=~~ \mathbf{z}_{j} - \Jac(\mathbf{z}_j)^+ \,
\phi(\mathbf{z}_j), ~~~~j = 0, 1, \ldots.
\end{equation}
where ~$\Jac(\mathbf{z})^+$ ~is the pseudo-inverse of the Jacobian
~$\Jac(\mathbf{z})$ ~of ~$\phi(\mathbf{z})$ ~given in
(\ref{eq:Jac_diffeo2}) ~and
~$\mathbf{z}_0 = (\alpha,\lb p_1 \rb, \ldots, \lb p_r \rb)$.

Detailed discussion on the Gauss-Newton iteration can be found in
\cite{Zeng09b,ZZ-tnt}.
~In a nutshell, the iteration (\ref{gni}) locally converges to the
least squares solution ~$\mathbf{z}_*$ ~that is the point satisfying
\begin{equation}  \label{phmin}
\big\|\phi(\mathbf{z}_*) - (\lb f\rb, 1, \ldots, 1)
\big\|_2 ~~=~~
\min_{ \mathbf{z} \in \C\times \C^{\langle \mathbf{m}_1\rangle} \times
\cdots \times \C^{\langle \mathbf{m}_r\rangle} }
 \big\|\phi(\mathbf{z}) - (\lb f\rb, 1, \ldots, 1)
\big\|_2
\end{equation}
if both the residual ~$\big\|\phi(\mathbf{z}_*) - (\lb f\rb,
1, \ldots, 1) \big\|_2$ ~and the initial error
~$\|\mathbf{z}_0 - \mathbf{z}_*\|_2$ ~are small.

\vspace{-3mm}
\begin{lemma}
~Let ~$\alpha f_1^{k_1}\cdots f_r^{k_r}$ ~be the numerical
squarefree irreducible factorization of ~$f$
~within ~$\epsilon$
~with ~$\mathbf{m}_j = \dg(f_j)$ ~for ~$j=1,\ldots,r$.
~Then, for almost all unit vectors
~$\mathbf{b}_i \in \C^{\langle \mathbf{m}_i \rangle}$,
~$i=1,\ldots,r$, ~there is a
factorization ~$\alpha_* f_{*1}^{k_1}\cdots f_{*r}^{k_r}\,\sim\,
\alpha f_1^{k_1}\cdots f_r^{k_r}$ ~such that
~$\mathbf{z}_* \,=\, \big( \alpha_*, \lb f_{*1} \rb, \ldots, \lb f_{*r}
\rb \big)$
~is the least squares solution to the equation {\em (\ref{nls})} with residual
\begin{equation} \label{reseq}
\big\|\phi(\mathbf{z}_*) - (\lb f\rb, 1, \ldots, 1)
\big\|_2 ~~=~~
\|f - \alpha f_1^{k_1}\cdots f_r^{k_r} \|
\end{equation}
where ~$\phi$ ~is defined as in {\em (\ref{eq:diffeo2})}.
\end{lemma}

\vspace{-3mm}
\proof
~Clearly ~$\big\|\phi(\mathbf{z}) - (\lb f\rb,1,\ldots,1)
\big\|_2  \ge
\|f - \alpha f_1^{k_1}\cdots f_r^{k_r} \|$
~by Definition~\ref{def:NIF}.
~On the other hand, setting ~$f_{*i} =
\frac{1}{\mathbf{b}_i \cdot \lb f_i \rb} f_i$
~for ~$i = 1, \cdots, r$ ~and an appropriate ~$\alpha_*$ ~yields
(\ref{reseq}).
\foorp

Under the main condition that the Jacobian ~$\Jac(\mathbf{z}_*)$ ~is
of full rank, the Gauss-Newton iteration converges locally \cite{Zeng09b}.
~The local convergence of the Gauss-Newton iteration requires two
conditions: ~The initial iterate ~$\mathbf{z}_0$ ~must be near the
least squares solution ~$\mathbf{z}_*$ ~and the residual ~$\|
\phi(\mathbf{z}_*) - (\lb f\rb, 1, \ldots, 1)\|_2$ ~must be sufficiently
small.
~From (\ref{reseq}), ~the residual ~$\|\phi(\mathbf{z}_*) -
(\lb f\rb, 1, \ldots, 1)\|_2$ ~is bounded by the data error
~$\|f-\tilde{f}\|$.
~As a result, ~the residual requirement will be satisfied if the data error
is sufficiently small.

The algorithmic and technical details of the numerical factorization are out
of the scope of this paper and will be elaborated in a separate works.

\vspace{-5mm}
\section{Implementation, software and sample results} \label{s:imp}
\vspace{-5mm}

Our numerical factorization algorithm is implemented for both
univariate and multivariate polynomials as a function {\tt
PolynomialFactor} in the Matlab package {\tt NAClab} for numerical
algebraic computation as an upgrade and an expansion from its
predecessor {\tt Apalab} \cite{Zeng-ApaTools}.
~The entire {\tt NAClab} package is freely available\footnote{\tt
http://homepages.neiu.edu/$\sim$naclab.html}, including numerical factorization,
numerical rank-revealing, numerical
computation of multiplicity structure at zeros of nonlinear systems,
numerical greatest common divisors, etc.
~We shall present several sample results highlighting the major improvement
areas of our algorithm and the resulting software: Efficiency, accuracy,
versatility and user friendliness.
~All the tests are carried out on a Samsung Series 7 XE700T1A tablet
computer with 4GB memory and Intel i5-2467M CPU at 1.60 GHz running on
Windows 7 64-bit operating system.
~The test log and relevant Matlab/Maple scripts can be downloaded
online\footnote{\tt http://homepages.neiu.edu/$\sim$zzeng/NumFactorTests.zip}.

The Matlab package {\tt NAClab} provides a user friendly
interface for numerical algebraic computations.
~Polynomials can be entered and output as intuitive strings for casual users.
~The function {\tt PolynomialFactor} can be conveniently executed as follows.

\scriptsize
\begin{verbatim}
    >> p = '-4 - 12*x*y + x^3*y^2*z + 3*x^4*y^3*z + 8*z^3 - 2*x^3*y^2*z^4'
    >> PolynomialFactor(p,1e-10,'row')
    ans =
    (-12) * (0.333333333333333 + x*y - 0.666666666666667*z^3) * (1 - 0.25*x^3*y^2*z)
\end{verbatim}
\normalsize

\vspace{-3mm}
\begin{example}[Univariate factorization]
\em
~Accurate factorization of univariate polynomials with multiple roots
has been a challenge in numerical computation.
~Conventional software functions for polynomial
root-finding, such as Matlab {\tt roots} and Maple {\tt fsolve}
can not factor such polynomial accurately and output scattered root clusters.
~For example, let
\begin{eqnarray} f(x) & = &
x^{100} - \mbox{\tiny $222.222222222222$}\,x^{99} +\cdots
- \mbox{\tiny $8.53544016536406\,10^{30}$}\,x +
\mbox{\tiny $1.47799829703274\,10^{29}$} \label{rt100} \\
& \approx &
(x-\mbox{\tiny $4.444444444444444$})^{10}
(x-\mbox{\tiny $3.333333333333333$})^{20}
(x-\mbox{\tiny $2.222222222222222$})^{30}
(x-\mbox{\tiny $1.111111111111111$})^{40} \nonumber
\end{eqnarray}

In contrast,  our {\tt PolynomialFactor} is an advanced polynomial
root-finder that is capable of accurate computation for multiple
roots without extending machine precision even if the coefficients
are perturbed.
~On this example, our {\tt PolynomialFactor}
yields a factorization containing accurate roots and multiplicities:

\scriptsize
\begin{verbatim}
    >> PolynomialFactor(f,1e-10,'row')
    ans =
    (x-4.44444444445)^10 * (x-3.33333333333)^20 * (x-2.22222222222)^30 * (x-1.11111111111)^40
\end{verbatim}
\normalsize
This is a substantial improvement over its predecessor \cite{zz05}.
\foorp
\end{example}
\vspace{-3mm}

Since available software implementations for multivariate
factorizations are built on different platforms, based on different
notions of numerical factorizations, and with different designing
emphases,
comprehensive comparisons are not feasible.
~Among them, Maple {\tt factor} is built for the exact factorization.
~Developed by Vershelde, the package {\tt PHC}\cite{PHC} is a general-purpose
polynomial system solver whose factorization option {\tt -f} is perhaps the
first implemented numerical factorization software.
~This {\tt PHC} option initiates the implementation of a factorization
algorithm \cite{SVW04} in numerical computation
based on the homotopy continuation method.
~The Maple code {\tt appfac} is developed by Kaltofen, May, Yang and Zhi
\cite{KMYZ08} and the algorithm uses similar reducibility test based on
\cite{Gao03,rup99}, which appears to be superior
in factoring polynomial with highly perturbed data.
~The computing examples in the remainder of this section are
designed to showcase the differences and improvement areas of our algorithm and
implementation.

\vspace{-3mm}
\begin{example}[Stewart-Gough Platforms]  \label{e:sg}
\em
~In \cite{SVW04}, the authors tested three polynomials derived
from the Stewart-Gough platform manipulator in mechanical engineering:
\begin{equation} \label{sgpolyn}
\left\{ \begin{array}{ccl}
g_1 & = & F_1(q_0,q_1,q_2,q_3)\,(q_0^2+q_1^2+q_2^2+q_3^2)^3 \\
g_2 & = & F_2(q_0,q_1,q_2,q_3)\,(q_0^2+q_1^2+q_2^2+q_3^2)^3 \\
g_3 & = & a p_3^3 (q_0 + bq_3)(q_0+cq_3)(q_0+i\,q_3)^5(q_0-i\,q_3)^5
\end{array} \right.
\end{equation}
where ~$F_1, \,F_2 \in \C[q_0,q_1,q_2,q_3]$.
~The polynomials ~$g_1$ ~and ~$g_2$ ~both have 910 terms while ~$g_3$ ~has 24.
~We test {\tt PHC} in windows 7 command prompt using the compiled executable
file {\tt phc.exe} provided by its authors compared with our interpretive
code {\tt PolynomialFactor} in Matlab.
~To level the base of accuracy comparison,
we disabled the Gauss-Newton iteration option in our {\tt PolynomialFactor}
in this test since the iterative refinement was not developed for the
computed factors when {\tt PHC}  was released.
~Table~\ref{tb:sg} lists the elapsed execution times and the errors,
where the forward errors are measured on the known factors only.
~Since the three implementations are tested on different platforms,
the comparisons should be considered indirect.
~Nonetheless, the results appears to show our algorithm is efficient
and accurate on those polynomials.
~It also appears that {\tt PHC} has been improved substantially as it runs
much faster and outputs more accurate factors than it is reported in 2004.
\end{example}

\begin{table}[ht]
\begin{center}
\begin{tabular}{||c|c||c|c|c||} \hline \hline
 \multicolumn{2}{||c||}{ } & $g_1$ & $g_2$ & $g_3$ \\ \hline \hline
\raisebox{-0mm}{Maple {\tt factor}}
 &       & \multicolumn{3}{c||}{not designed for empirical data} \\
\hline
{\tt appfac}
 & \multicolumn{1}{c||}{ }  & \multicolumn{3}{c||}{$---------$}
\\ \hline
& elapsed time
&\scriptsize 1382.8
&\scriptsize 1410.1
&\scriptsize 1.48  \\
{\tt PHC} & backward error
&\scriptsize $9.6\times 10^{-10}$
&\scriptsize 0.9945
&\scriptsize $1.7\times 10^{-12}$ \\
 & forward error
&\scriptsize $1.3 \times 10^{-12}$
&\scriptsize $1.1\times 10^{-12}$
& \scriptsize $3.4\times 10^{-13}$ \\
\hline
&\bf elapsed time
&\scriptsize\bf 376.5
&\scriptsize\bf 480.3
&\scriptsize\bf 0.79 \\
{\footnotesize \tt \bf PolynomialFactor}
 &\bf backward error
&\scriptsize $\mathbf{4.3\times 10^{-15}}$
&\scriptsize $\mathbf{4.1\times 10^{-15}}$
&\scriptsize $\mathbf{4.8\times 10^{-14}}$ \\
\scriptsize (without refinement)
&\bf forward error
&\scriptsize $\mathbf{2.0\times 10^{-15}}$
&\scriptsize $\mathbf{8.5\times 10^{-16}}$
&\scriptsize $\mathbf{5.0\times 10^{-16}}$  \\
\hline \hline
\end{tabular}
\end{center}
   \vspace{-5mm}
\caption{\small Factorization results on polynomials in (\ref{sgpolyn})
derived from Stewart-Gough platforms.}
\label{tb:sg}
\end{table}

\vspace{-3mm}
\vspace{-3mm}
\begin{example}[A polynomial with 5 numerical factorizations]
\label{ex:5f}
\em
~An issue of significant importance on the concept of numerical
factorization is that
a polynomial may have different numerical factorizations within different
error tolerances approximating different conventional factorizations.
~A numerical factorization algorithm in this context needs
mechanisms for targeting specific factorizations.
~For example, the polynomial

\vspace{-8mm}
\footnotesize
\begin{eqnarray}
p_1 &~=~& x^7\,y +x^5\,y^3 +x^6 -x\,y^7 -x^3\,y^5
-\mbox{\scriptsize \tt 3}\,x^2\,y^4
+\mbox{\scriptsize \tt 6}\,x^2\,y^2
+2(x^3\,y^3 -x^2\,y^6 +x^6\,y^2 -x\,y^5-y^6 +x^4) \nonumber \\
& &
+\mbox{\scriptsize \tt 7}\,x^4\,y^2
+\mbox{\scriptsize \tt 4}\,x^5\,y
+\mbox{\scriptsize \tt .999001}\,x\,y^3
+\mbox{\scriptsize \tt 1.998002001}\,x\,y
+\mbox{\scriptsize \tt 4.999001}\,x^3\,y
+\mbox{\scriptsize \tt 2.999001}\,y^2
-\mbox{\scriptsize \tt 1.000999}\,x^2
\nonumber \\
& &
+\mbox{\scriptsize \tt .001}(y^3+x^3
+\mbox{\scriptsize \tt 3}\,x\,y^2
+x^3\,y^2
+\mbox{\scriptsize \tt 3}\,x^2\,y +x^4\,y +x\,y^4 +x^2\,y^3
+\mbox{\scriptsize \tt 2}\,y
+\mbox{\scriptsize \tt 2}\,x)
-\mbox{\scriptsize \tt 2.001997998999} \label{p1}
\end{eqnarray}
\normalsize
can be considered as empirical data of ~$p_1$ ~itself and
any one of the four factorable polynomials
\footnotesize
\begin{equation}  \label{p2345}
\begin{array}{ccl}
p_2 & = &
(xy+\mbox{\scriptsize \tt 1})\,(
\mbox{\scriptsize \tt  .001} x^2 y+\mbox{\scriptsize \tt 6} x^3 y+
\mbox{\scriptsize \tt 2} x^5 y+\mbox{\scriptsize \tt  .001} x y^2
+\mbox{\scriptsize \tt 2} x^2 y^2+\mbox{\scriptsize \tt  .001} y^3
-\mbox{\scriptsize \tt 1.000999} x^2+\mbox{\scriptsize \tt 2} x^4 - y^6+x^6\\
& & ~~+x^4 y^2+ \mbox{\scriptsize \tt 4} x y-x^2 y^4+\mbox{\scriptsize \tt .001} x^3
+\mbox{\scriptsize \tt 2.999001} y^2-\mbox{\scriptsize \tt 2} x y^5
-\mbox{\scriptsize \tt 2.001997999}+ \mbox{\scriptsize \tt .002}  y+
\mbox{\scriptsize \tt .002} x-\mbox{\scriptsize \tt 2} x y^3) \\
p_3 & = & (-\mbox{\scriptsize \tt 2} x y^3+\mbox{\scriptsize \tt 2} x y
+\mbox{\scriptsize \tt 2} x^3 y+x^4+\mbox{\scriptsize \tt 2} y^2-y^4
-\mbox{\scriptsize \tt 1.000999}+\mbox{\scriptsize \tt .001} y
+\mbox{\scriptsize \tt .001} x
)
(x^2+y^2+\mbox{\scriptsize \tt 2})(xy+\mbox{\scriptsize \tt 1})
\\
p_4 & = & (x^3-x y^2+x+x^2 y-y^3+y+x^2-y^2+ \mbox{\scriptsize \tt 1.001})(x+y-\mbox{\scriptsize \tt 1})
(x^2+y^2+\mbox{\scriptsize \tt 2})(xy+\mbox{\scriptsize \tt 1})
\\
p_5 & = & (x+y+\mbox{\scriptsize \tt 1})(x^2-y^2+\mbox{\scriptsize \tt 1})
(x+y-\mbox{\scriptsize \tt 1})
(x^2+y^2+\mbox{\scriptsize \tt 2})(xy+\mbox{\scriptsize \tt 1})
\end{array}
\end{equation} \normalsize
with data errors of various magnitudes listed in Table~\ref{tb:p1-5}.
~In other words, the polynomial ~$p_1$ ~has a numerical factorization
~$1\cdot p_1$ ~within an error tolerance between ~$0$ ~and ~$10^{-14}$,
~and numerical factorizations within error tolerances roughly in the
intervals ~$(10^{-13},10^{-12})$, ~$(10^{-9},10^{-8})$, ~$(10^{-6}, 10^{-5})$
~and ~$(10^{-3},10^{-2})$  ~approximating the exact factorizations
of ~$p_2$, ~$p_3$, ~$p_4$, ~$p_5$ ~in (\ref{p2345}) respectively.
~Our formulation of the numerical factorization includes the error tolerance
~$\epsilon$ and our implementation {\tt PolynomialFactor} provides such an
option.
~Based on the choices of those error tolerances,
{\tt PolynomialFactor} calculates all five
numerical factorizations (\ref{p2345})
with forward accuracies in the same orders of the data errors as shown
in Table~\ref{tb:p1-5}.
~Other algorithms such as Maple {\tt factor}, {\tt PHC} and {\tt appfac} are
designed to compute one numerical factorization from a given polynomial data.

\end{example}

\begin{table}[ht]
\begin{center}
\begin{tabular}{||l||c|c|c|c|c||} \hline \hline
underlying polynomials & & & & & \\
~~to be factored & \raisebox{2mm}{$p_1$} & \raisebox{2mm}{$p_2$}
& \raisebox{2mm}{$p_3$} & \raisebox{2mm}{$p_4$} & \raisebox{2mm}{$p_5$}
\\ \cline{1-1}
data error ~$\sin(p_1,p_j)$ &\scriptsize  $0$
&\scriptsize  $7.36\times 10^{-14}$
&\scriptsize  $1.05\times 10^{-10}$
&\scriptsize  $2.56\times 10^{-7}$
&\scriptsize  $4.92\times 10^{-4}$  \\
\hline \hline
\footnotesize Maple {\tt factor}
& \scriptsize $0$
& \scriptsize $----$
& \scriptsize $----$
& \scriptsize $----$
& \scriptsize $----$
\\ \hline
{\footnotesize {\tt appfac}+refinement}
& \scriptsize $----$
& \scriptsize $----$
& \scriptsize $----$
& \scriptsize $2.38\times 10^{-7}$
& \scriptsize $----$
\\ \hline
{\footnotesize {\tt PHC}} {\scriptsize (with no refinement)}
& \scriptsize $----$
& \scriptsize $----$
& \scriptsize $1.43\times 10^{-9}$
& \scriptsize $----$
& \scriptsize $----$
\\ \hline
{\footnotesize \tt \bf PolynomialFactor}
& \scriptsize $\mathbf{6.96\times 10^{-16}}$
& \scriptsize $\mathbf{8.88\times 10^{-14}}$
& \scriptsize $\mathbf{9.16\times 10^{-11}}$
& \scriptsize $\mathbf{2.06\times 10^{-7}}$
& \scriptsize $\mathbf{5.62\times 10^{-4}}$ \\ \hline \hline
\end{tabular}
\end{center}
   \vspace{-5mm}
\caption{\small Forward accuracies of
Maple {\tt factor}, {\tt PHC}, {\tt appfac} and
{\tt PolynomialFactor} from the data polynomial ~$p_1$ ~in (\ref{p1})
calculating the factorization of either ~$p_1,p_2,p_3,p_4$ ~or ~$p_5$
(\ref{p2345}).}
\label{tb:p1-5}
\end{table}

\vspace{-5mm}
\section{Conclusions}
\label{s:conc}

\vspace{-5mm}
Conventional factorization is an ill-posed problem in the sense that it is
infinitely sensitive to data perturbations.
~The reason for such hypersensitivity is revealed by the geometry of polynomial
factorization and the singularity can be quantified by the dimension deficit of
the factorization manifold.
~The numerical factorization as formulated in this paper generalizes the
concept of conventional factorization and eliminates the ill-posedness.
~By establishing the fundamental theorems for geometric structure of
multivariate factorization, we proved that the numerical
factorization uniquely exists and possesses Lipschitz continuity with respect
to data under the overall assumption that the data error is small.
~Consequently, the numerical factorization achieves the objective of
recovering the factorization accurately even if the polynomial data are
empirical and the accuracy is in the order of data precision.
~An algorithm is implemented as a Matlab module and numerical results support
this conclusion.

\bibliographystyle{plain}

\end{document}